\documentclass[11pt]{article}
\usepackage{amsmath}
\usepackage{graphicx}
\usepackage{epsfig}
\usepackage{amsfonts}
\usepackage{amssymb}
\usepackage{placeins}
\usepackage{cases}
\usepackage[margin=1in]{geometry}
\usepackage{authblk}
\usepackage[titletoc,toc,title]{appendix}
\usepackage{epstopdf}

\usepackage{accents}

\title{Optimal Control of Oscillation Timing and Entrainment Using Large Magnitude Inputs:~An Adaptive Phase-Amplitude-Coordinate-Based Approach}

\begin{document}
\author[1]{Dan Wilson}
\affil[1]{Department of Electrical Engineering and Computer Science, University of Tennessee, Knoxville, TN 37996, USA}
\maketitle


\begin{abstract}
Given the high dimensionality and underlying complexity of many oscillatory dynamical systems, phase reduction is often an imperative first step in control applications where oscillation timing and entrainment are of interest.  Unfortunately, most phase reduction frameworks place restrictive limitations on the magnitude of allowable inputs, limiting the practical utility of the resulting phase reduced models in many situations.  In this work, motivated by the search for control strategies to hasten recovery from jet-lag caused by rapid travel through multiple time zones, the efficacy of the recently developed adaptive phase-amplitude reduction is considered for manipulating oscillation timing in the presence of a large magnitude entraining stimulus.  The adaptive phase-amplitude reduced equations allow for a numerically tractable optimal control formulation and the associated optimal stimuli significantly outperform those resulting from from previously proposed optimal control formulations.  Additionally, a data-driven technique to identify the necessary terms of the adaptive phase-amplitude reduction is proposed and validated using a model describing the aggregate oscillations of a large population of coupled limit cycle oscillators.   Such data-driven model reduction algorithms are essential in situations where the underlying model equations are either unreliable or unavailable.  
\end{abstract}

\section{Introduction}
Mathematical analysis and control of oscillatory dynamical systems is a widely studied problem with relevant applications to neurological brain rhythms  \cite{nabi132}, \cite{holt16}, \cite{kope14}, circadian physiology \cite{diek16}, \cite{bagh08}, \cite{serk14}, and various other physical and chemical systems \cite{stro05}, \cite{dorf13}, \cite{dell20}, \cite{wils18bz}, \cite{zlot16}.  Given the high dimensionality and sheer complexity of many oscillatory dynamical systems, model reduction is often a necessary first step in their analysis.  While many reduction strategies have been developed for limit cycle oscillators, phase reduction \cite{winf01}, \cite{erme10}, \cite{kura84} is one of the most widely applied techniques which allows for oscillatory systems of the form
\begin{equation} \label{basicphase}
\dot{x} = F(x) + U(t),
\end{equation}
where $x \in \mathbb{R}^N$ is the system state, $F$ describes the nominal dynamics, and $U$ is an exogenous input to be analyzed according to
\begin{equation} \label{predintro}
\dot{\theta} = \omega + Z(\theta) \cdot U(t),
\end{equation}
where the phase $\theta \in [0,2\pi)$ characterizes the state of the oscillator in reference to a limit cycle, $\omega(p_0) = 2\pi/T$ is the natural frequency  where $T$ is the nominal, unperturbed period, $Z(\theta)$ is the phase response curve (i.e.,~the gradient of the phase coordinate evaluated on the limit cycle), and the `dot' denotes the dot product.   Phase reduction provides a tremendous decrease in the dimensionality of a given limit cycle oscillator, allowing for the original $N$-dimensional, generally nonlinear behavior to be studied as a 1-dimensional ordinary differential equation.  In exchange for this tremendous reduction in dimensionality, the magnitude of the perturbations is required to be uniformly bounded in time by $\epsilon$, where $0 < \epsilon \ll 1$ to ensure that Equation \eqref{predintro} is valid to first order accuracy in $\epsilon$. Practically, the reduction \eqref{predintro} can be applied to make predictions about system behavior in response to larger magnitude inputs with the understanding that its performance will begin to degrade as the magnitude of the input becomes larger; precise bounds on the allowable $U(t)$ are usually not known {\it a priori} and vary on an application-by-application basis.

In an effort to extend the applicability of phase reduction when larger magnitude inputs are considered, a variety of reduction algorithms have been proposed that take into account the transient dynamics in directions transverse to the periodic orbit.  For example, the notions of entrainment maps \cite{diek16}, operational phase coordinates \cite{wils18operat}, functional phase response curves \cite{cui09}, local orthogonal rectification \cite{lets20}, stochastic phase \cite{bres18}, and average isophase \cite{schw13} rely on different definitions of `phase', with each being well-suited for specific applications.  Other recent works have attempted to use phase-amplitude coordinate systems to identify  reduced order equations that are valid beyond first order accuracy in $\epsilon$, \cite{wils17isored}, \cite{geng20}, \cite{wils20highacc}, \cite{rose19}, \cite{wils19prl}.  While such techniques generally provide better results than comparable first order accurate reduction methods, they still suffer from limitations that preclude the consideration of medium-to-large magnitude inputs.  Koopman analysis is an emerging reduction framework that can be used to represent a nonlinear dynamical oscillator (or more generally any nonlinear dynamical system) with a linear, but infinite dimensional operator \cite{budi12}, \cite{mezi19}, \cite{kais17}. However, it can be difficult to identify a suitable finite dimensional basis that fully captures the dynamics of the Koopman operator.

Recently, an adaptive phase-amplitude reduction strategy was proposed in \cite{wils21adaptive} that uses the standard definition of asymptotic phase based on isochrons \cite{winf01}, \cite{guck75} and also considers the slowly decaying isostable coordinates \cite{wils16isos} which represent level sets of the slowest decaying Koopman eigenfunctions \cite{maur13}.  By considering a family of limit cycles that emerge for different parameter sets, the nominal parameter set can be adaptively chosen in order to limit the error in the reduced order equations.  As shown in \cite{wils21adaptive}, provided the Floquet exponents of the truncated isostable coordinates are $O(1/\epsilon)$ terms, the adaptive phase-amplitude reduction is valid to $O(\epsilon)$ accuracy even when the input $U(t)$ is an $O(1)$ term.  

Because the adaptive phase-amplitude reduction framework allows for the consideration of particularly large magnitude inputs, this reduction framework is an attractive candidate for control of oscillation timing in situations where the application of large magnitude inputs is unavoidable.  In this work, optimal control of oscillation timing in the presence of a large magnitude entraining stimulus will be considered.  This situation can be used to represent circadian rhythms that are entrained to a 24-hour light-dark cycle.  In humans, a population of roughly 10,000 coupled neurons referred to as the suprachiasmatic nucleus is responsible for maintaining circadian time \cite{moor02}, \cite{repp02}.  The body's circadian cycle is nominally entrained to a 24-hour light-dark cycle in order to regulate sleep and wake timing among other physiological processes.  Circadian misalignment, more commonly referred to as jet-lag, represents a disruption to this nominal steady entrainment and results from a mismatch between the environmental time and one's own circadian clock \cite{sack10}.   

Numerical models of circadian oscillations have been helpful for developing post-travel treatment protocols to aid in recovery from circadian misalignment \cite{bagh08}, \cite{serk14}, \cite{dean09}, however, optimal control algorithms can typically only be implemented on in simplified, low-dimensional models.  For more complicated and high-dimensional models, standard phase reduction techniques typically fail.  There are two fundamental reasons for this:~the first is that circadian rhythms result from the collective oscillation of a large number of coupled oscillators.  Phase reduction methods can be used for these population oscillations using techniques described in \cite{ko09}, \cite{leva10}, \cite{kawa08} and \cite{dumo17}, however, the inter-oscillator coupling usually results in amplitude coordinates that decay slowly.  Practically, these amplitude coordinates represent shifts in the distribution of phases of the individual oscillators;  without the ability to accurately account for these time-varying phase distributions only particularly small magnitude inputs can be considered.  The second difficulty that arises when using phase reduction methods for applications involving circadian oscillations is that the nominal 24-hour light-dark cycle itself is generally be considered a strong perturbation  that shifts the oscillator state far from its nominal limit cycle \cite{diek18}.   For all but a small subset of oscillator models with dynamics that can be greatly simplified when certain assumptions are made \cite{kura84}, \cite{ott08},  it is difficult to predict how the aggregate behavior will respond to the large magnitude entraining stimuli precluding the consideration of control inputs to manipulate the phase.  For these reasons, reduction techniques that allow for large inputs are essential for applications that involve circadian oscillations.  

There are two primary focal points of this work.  The first is the development and assessment of strategies to implement the adaptive phase-amplitude reduction in an optimal control framework to manipulate oscillation timing and speed recovery from circadian misalignment.  In the results to follow, by making appropriate assumptions about the characteristics of the adaptive reduction, a calculus of variations approach can be used to identify optimal control inputs in the presence of a large magnitude entraining stimulus.  This strategy can be used to identify inputs that engender significantly larger time shifts than other recently proposed oscillation timing control frameworks.  The secondary focus of this work is the development of data-driven techniques for identification of the necessary terms of the adaptive phase-amplitude reduction that are valid in situations where numerical models are either unreliable or unavailable.  By leveraging techniques developed in \cite{wils20ddred} a general strategy is proposed to accurately determine necessary terms of the adaptive reduction solely from observable data.   

The organization of this paper is as follows:~Section \ref{backsec} gives necessary background information about the adaptive phase-amplitude reduction method used to identify reduced order models in this work.  Section \ref{optsec} considers an optimal control formulation using the adaptive phase-amplitude reduction for  manipulation of oscillation timing in the presence of an external entraining stimulus.  Analogous, previously proposed formulations are also considered that use phase-only and first order accurate phase-amplitude reductions in order to provide comparisons with the proposed optimal control formulation.  Section  \ref{simpleresult} provides results for the optimal control framework when using a relatively simple, three-dimensional model of  circadian oscillations.  Section \ref{popsec} considers a population-level model of coupled oscillators that give rise to a collective oscillation.  Here, a strategy is presented for identifying the necessary terms of the adaptive reduction from a single system observable and the optimal control strategy is subsequently implement to identify efficient stimuli to hasten recovery from circadian misalignment.  Section \ref{concsec} provides concluding remarks.

\section{Background} \label{backsec}
In the analysis to follow, a phase reduction \eqref{predintro} will be considered that explicitly incorporates a set of nominal parameters.  In this context, consider an ordinary differential equation of the form
\begin{equation}\label{maineq}
\dot{x} = F(x,p_0) + U(t),
\end{equation}
where $x$, $F$, and $U$ are defined identically to the terms from \eqref{basicphase} and $p_0 \in \mathbb{R}^M$ is a set of nominal parameters.   Suppose that for a constant choice of $p_0$, \eqref{maineq} has a stable periodic orbit $x^\gamma_{p_0}$.   Provided that $U(t)$ is sufficiently small,  phase reduction is a well established technique \cite{erme10}, \cite{kura84}, \cite{winf01} that can be used to analyze the behavior of \eqref{maineq} in the weakly perturbed limit according to
\begin{equation} \label{standardphase}
\dot{\theta} = \omega(p_0) + Z(\theta,p_0) \cdot U(t),
\end{equation}
where $\omega$ and $Z$ are defined identically to the terms from \eqref{predintro}.

\subsection{Isostable Coordinate Reduced Frameworks}

In situations where larger magnitude inputs must be considered, phase reduction can still be used but one must also consider the transient dynamics in directions transverse to the limit cycle.  Various phase-amplitude coordinate frameworks have been developed for this purpose \cite{wils17isored}, \cite{shir17}, \cite{wedg13}, \cite{lets20}, \cite{bres18}, \cite{wils18operat}, \cite{rose19}; this work will focus on strategies that use isostable coordinates to characterize the amplitude coordinates \cite{maur13}, \cite{wils17isored}, which represent level sets of the slowest decaying eigenfunctions of the Koopman operator \cite{maur13}.  Each isostable coordinate represents the magnitude of a particular Koopman eigenmode.  Intuitively, states that correspond to larger magnitude isostable coordinates will take longer to decay to the periodic orbit.  For the slowest decaying Koopman eigenfunctions, explicit definitions can be given for the isostable coordinates in the basin of attraction of a fixed point \cite{maur13} or periodic orbit \cite{wils17isored}, however, faster decaying isostable coordinates must be defined implicitly according to their respective Koopman eigenfunctions with decay rates that are governed by their associated Floquet exponents.  Previous work \cite{wils16isos} (cf.~\cite{cast13}) investigated strategies for computing the behavior of the isostable coordinates for models of the form \eqref{maineq} in response to perturbations in a neighborhood of the limit cycle yielding reduced order equations of the form
\begin{align} \label{isored}
\dot{\theta} &= \omega(p_0) + Z(\theta,p_0) \cdot U(t), \nonumber \\
\dot{\psi_j} &= \kappa_j(p_0) \psi_j + I_j(\theta,p_0) \cdot U(t), \nonumber \\
x(\theta,\psi_1,\dots,\psi_\beta) &= x^\gamma_{p_0}(\theta) + \sum_{k = 1}^\beta g^k(\theta) \psi_j, \nonumber \\
j &= 1,\dots,\beta,
\end{align}
where $\psi_j$ is the $j^{\rm th}$ isostable coordinate with corresponding Floquet exponent $\kappa_j$, $g^k(\theta)$ are associated Floquet eigenfunctions, and $I_j(\theta,p_0)$ is the isostable response curve that represents the gradient of the isostable coordinates evaluated on the periodic orbit.  In \eqref{isored}, the rapidly decaying amplitude coordinates are generally truncated so that only $\beta \leq N-1$ of the slowest decaying isostable coordinates are explicitly considered -- isostables with rapid exponential decay can generally be assumed to be zero without adverse effects in the accuracy of the reduction.  Numerical computation of $Z$, and each $I_j$ and $g^k$ can be performed using the `adjoint method' \cite{brow03} and related equations described in \cite{wils20highacc}.  The so-called `direct method' \cite{izhi07}, \cite{neto12} and related data-driven techniques \cite{wils19phase}, \cite{wils20ddred} have been developed for inferring the necessary terms of \eqref{isored} from experimental data when the underlying model equations are unknown.

Much like \eqref{standardphase}, Equation \eqref{isored} is only valid provided the state remains sufficiently close to the periodic orbit.  Unlike \eqref{standardphase} alone, however, \eqref{isored} can be used to provide information about the transient behaviors that characterize the amplitude dynamics.  It is generally assumed that under the application of  input with $O(\epsilon)$ magnitude that each $\psi_j$ remains an $O(\epsilon)$.  In this case, \eqref{isored} is valid to leading order $\epsilon$.  Recent work has investigated isostable reduced equations that are valid to second order accuracy \cite{wils17isored}, \cite{wils19phase} and beyond \cite{wils20highacc}, however, these reduction frameworks still require the magnitude of the applied inputs to be sufficiently small.  Other reduction frameworks have been developed that are valid for inputs with arbitrary magnitude provided that they are sufficiently slowly varying \cite{kure13}, \cite{park16} or rapidly varying \cite{pyra15}, \cite{vela03}.  Nevertheless, few general reduction frameworks are available that are valid for arbitrary, large magnitude inputs. 

\subsection{Adaptive Phase-Isostable Reduction}
Isostable reduced equations of the form \eqref{isored} typically assume that the nominal system parameters are constant in the analysis of a perturbed limit cycle. Recent work \cite{wils21adaptive} has leveraged the phase-isostable coordinate reduction framework \eqref{isored} to develop an adaptive reduction that actively sets the nominal system parameters in an effort to keep isostable coordinates low, thereby allowing for inputs without $O(\epsilon)$ constraints (either in magnitude or rate of change).    As discussed in \cite{wils21adaptive}, to implement an adaptive phase-isostable reduction, first suppose that in some allowable range of nominal system parameters $p \in \mathbb{R}^M$ that $x_p^\gamma(t)$ is continuously differentiable with respect to both $t$ and $p$.  One can then rewrite Equation \eqref{maineq} as 
\begin{equation} \label{extendedmain}
\dot{x} = F(x,p) + U_e(t,p,x),
\end{equation}
where 
\begin{equation} \label{uedef}
U_e(t,p,x) = U(t) + F(x,p_0) - F(x,p).
\end{equation}

Intuitively, $F(x,p)$ represents the underlying dynamics for a given choice of system parameters  $p$, and $U_e$ represents the effective input.  Allowing $p$ to be variable, and assuming that $\theta(x,p)$ and $\psi_j(x,p)$ are continuously differentiable in a neighborhood of the limit cycle for all $x$ and $p$, one can transform \eqref{extendedmain} to phase and isostable coordinates via the chain rule
\begin{align} \label{fulleq}
\frac{d \theta}{dt}  &= \frac{\partial \theta}{\partial x} \cdot  \frac{d x}{dt} + \frac{\partial \theta}{\partial p} \cdot   \frac{d p}{dt}, \nonumber \\
\frac{d \psi_j}{dt } &= \frac{\partial \psi_j}{\partial x} \cdot  \frac{d x}{dt} + \frac{\partial \psi_j }{\partial p} \cdot  \frac{d p}{dt},
\end{align}
for $j = 1,\dots,\beta$.  Above, provided that the state is close to the periodic orbit $x^\gamma_p$,  $ \frac{\partial \theta}{\partial x} \cdot  \frac{d x}{dt} = \omega(p) + Z(\theta,p) \cdot U(t)$ and  $\frac{\partial \psi_j}{\partial x} \cdot  \frac{d x}{dt} = \kappa_j(p) \psi_j + I_j(\theta,p) \cdot U(t)$  as given by Equation \eqref{isored}.  For the remaining terms, as explained in \cite{wils21adaptive}, one can write
\begin{align} \label{dandq}
D(\theta,p) \equiv \frac{\partial \theta}{\partial p} = \begin{bmatrix}  Z(\theta,p) \cdot \frac{\partial \Delta x_p}{\partial p_1} & \dots &  Z(\theta,p) \cdot \frac{\partial \Delta x_p}{\partial p_M}  \end{bmatrix}^T, \nonumber \\
Q_j(\theta,p) \equiv \frac{\partial \psi_j}{\partial p} = \begin{bmatrix}  I_j(\theta,p) \cdot \frac{\partial \Delta x_p}{\partial p_1} & \dots &  I_j(\theta,p) \cdot \frac{\partial \Delta x_p}{\partial p_M}  \end{bmatrix}^T,
\end{align}
where 
\begin{equation} \label{partialorb}
\frac{\partial  \Delta x_p}{\partial p_i} = \lim_{dp_i \rightarrow 0}  \frac{x^\gamma_{p} - x^\gamma_{p+dp_i}}{dp_i}.
\end{equation}
Note that Equation \eqref{partialorb} represents the change in the nominal periodic orbit in response to a change in the parameter $p_i$.   In Equation \eqref{dandq} and \eqref{partialorb}, all derivatives are evaluated at $p$ and $\theta$ on the limit cycle $x^\gamma_p$.  Substituting both \eqref{dandq} and the phase and isostable dynamics from \eqref{isored} into \eqref{fulleq} yields the adaptive phase-isostable reduction
\begin{align} \label{fulladaptive}
\dot{\theta} &= \omega(p) +Z(\theta,p)   \cdot U_e(t,p,x)   +   D(\theta,p)     \cdot \dot{p},   \nonumber \\
\dot{\psi}_j &= \kappa_j(p) \psi_j +  I_j(\theta,p)  \cdot U_e(t,p,x)  +   Q_j(\theta,p)  \cdot \dot{p},  \nonumber \\
j &= 1, \dots, \beta, \nonumber \\
\dot{p} &= G_p(p,\theta,\psi_1,\dots, \psi_\beta),
\end{align}
where the function $G_p$ is designed to actively choose $p$ (and by extension $\dot{p}$) in a manner that keeps the isostable coordinates small.  As explained in \cite{wils21adaptive}, provided $G_p$ can be designed so that $\psi_j$ remain $O(\epsilon)$ for $j \leq \beta$ and that the neglected isostable coordinates have sufficiently large magnitude Floquet exponents, Equation \eqref{fulladaptive} is accurate to $O(\epsilon)$ provided that $U_e$ is an $O(1)$ term.  Recalling that the standard phase-isostable reduction \eqref{isored} assumes the input is an $O(\epsilon)$ term, the adaptive isostable reduction represents a substantial improvement that allows for significantly larger magnitude inputs to be considered.  It is worth emphasizing that the underlying model \eqref{maineq} does not need to have time-varying parameters in order for the adaptive reduction \eqref{fulladaptive} to be implemented.  Rather, the adaptive reduction actively changes the nominal parameter set $p$  so that the state is close to the periodic orbit $x^\gamma_{p}$ thereby limiting the magnitude of the isostable coordinates.

\section{Optimal Control of Oscillation Timing using Adaptive Phase-Amplitude Reduction} \label{optsec}

Both phase reduction \eqref{standardphase} and phase-amplitude reduction have been applied fruitfully to applications that involve control and analysis of behaviors that emerge in limit cycle oscillators \cite{stro05}, \cite{abra04}, \cite{wils19prl}, \cite{zlot16}, \cite{piet19}.  The effective reduction in dimension that results from phase reduction can allow for the application of control and analysis techniques that would otherwise be intractable.  However, the results are only  applicable when sufficiently weak perturbations are applied limiting the practical utility in many applications.  Here, a prototype problem is considered for identifying an optimal control input to advance or delay the oscillation timing in an oscillator subject to an external entraining stimulus.  This problem is motivated by the search for jet-lag mitigation protocols that can be used to limit the negative effects of circadian misalignment by allowing one's circadian cycle to reentrain rapidly to a new time zone \cite{bagh08}, \cite{serk14}, \cite{dean09}.   Related problems were considered in \cite{moeh06} and \cite{nabi132} using standard phase reduction methods, and in \cite{mong19b} (resp.~\cite{wils17isored}) using first (resp.~second)  order accurate phase-amplitude reduced equations.  As shown in the results to follow, the problem formulation using the adaptive phase-amplitude reduced equations allows for inputs that are far larger in magnitude than those considered previously without sacraficing accuracy -- consequently, inputs that provide significantly larger magnitude changes in oscillation timing can be obtained.

To begin, consider the adaptive phase-isostable reduced equations from \eqref{fulladaptive} with only a single isostable coordinate
\begin{align} \label{adaptiveone}
\dot{\theta} &= \omega(p) +Z(\theta,p)   \cdot U_e(t,p,x)   +   D(\theta,p)     \dot{p},   \nonumber \\
\dot{\psi} &= \kappa(p) \psi +  I(\theta,p)  \cdot U_e(t,p,x)  +   Q(\theta,p)  \dot{p},  \nonumber \\
\dot{p} &= G_p(p,\theta,\psi).
\end{align}
Above, for convenience of notation, the subscripts on the terms involving isostable coordinates have been omitted because there is only one isostable coordinate.  Additionally, it will be assumed that $p \in \mathbb{R}^1$.  In the analysis to follow, it will be assumed that the input can be written as $U(t) = \delta u(t)$, where $\delta \in \mathbb{R}^N$ is a constant vector and $u(t) \in \mathbb{R}$ is a time-varying control signal.  In other words,  $U(t)$ a rank-1 perturbation.  It will also be assumed that $F(x,p) = F(x) + p \delta$.  This situation results, for instance, when allowing the input to also be considered as a time-varying parameter in the adaptive reduction.  Under these assumptions, Equation \eqref{adaptiveone} simplifies to
\begin{align} \label{adaptivemod}
\dot{\theta} &= \omega(p) + z(\theta,p)  (u(t) + p_0 - p)   +   D(\theta,p)     \dot{p},   \nonumber \\
\dot{\psi} &= \kappa(p) \psi +  i (\theta,p)  (u(t) + p_0 - p)  +   Q(\theta,p)  \dot{p},  \nonumber \\
\dot{p} &= G_p(p,\theta,\psi),
\end{align}
where $z(\theta,p) = Z(\theta,p) \cdot \delta$ and $i(\theta,p) = I(\theta,p) \cdot \delta$.  To further simplify \eqref{adaptivemod}, it will be explicitly assumed that $|Q(\theta,p)| > \nu$ for all allowable $\theta$ and $p$ with $\nu>0$.  In this situation, taking $G_p(p,\theta,\psi) = -R(\theta,p) (u  + p_0 - p)$ where $R(\theta,p) = \frac{i(\theta,p)}{Q(\theta,p)}$, the isostable dynamics from equation \eqref{adaptivemod} simplify to $\dot{\psi} = \kappa(p) \psi$ which has a single globally stable equilibrium at $\psi = 0$.  Noting that $G_p$ does not depend on $\psi$, the isostable coordinate dynamics can be ignored from \eqref{adaptivemod} yielding
\begin{align} \label{adaptiveframe}
\dot{\theta} &= \omega(p) + z(\theta,p)  (u(t)  - p)   -   D(\theta,p)   R(\theta,p) (u   - p),   \nonumber \\
\dot{p} &= -R(\theta,p) (u  - p),
\end{align} 
where $p_0 = 0$ is assumed for simplicity of exposition.  To formulate the optimal control problem, consider a $T_{\rm ext}$-periodic external input $u(t) = u_{\rm nom}(t_s)$ applied to Equation \eqref{maineq} where
\begin{equation} t_s = 
\begin{cases}
t, &\mbox{for } t \leq t_0, \\
t + \Delta t, & \mbox{for } t > t_0.
\end{cases}
\end{equation}
Suppose that for $t\leq t_0$, the oscillator is fully entrained to the periodic orbit so that $x(t) = x(t + T_{\rm ext})$ for all $t + T_{\rm ext} < 0$.  This fully entrained orbit will be denoted by $x^\gamma_{{\rm ent}}(t_s)$.   In the context of circadian oscillations, $t_s$ would represent the environmental time that controls a 24-hour light-dark cycle $u_{\rm nom}(t_s)$.  Suppose that at $t = t_0$,  the environmental time instantaneously shifts by $\Delta t$ time units, for instance, representing a sudden shift across multiple time zones.   The control objective is to identify a control input $u(t) = u_{\rm nom}(t + \Delta t) + \Delta u(t)$ so that $x(t_0 + T_f) = x^\gamma_{p_0,{\rm ent}} (t_0 + \Delta t)$ that minimizes the cost functional  $C = \int_{t_0}^{t_0 + T_f}  \Delta u^2(t)  dt$ subject to the constraints on the allowable control $\Delta u_{\rm min} \leq  u(t) + u_{\rm nom}(t + \Delta t) \leq \Delta u_{\rm max}$.  Intuitively, this control problem seeks to identify an external input $\Delta u(t)$ so that the system is fully entrained after $T_f$ time units to the time-shifted entraining stimulus.

Working in the adaptive reduction framework from \eqref{adaptiveframe} this control problem can be approached by first defining the Hamiltonian associated with the cost functional
\begin{align} \label{hameq}
\mathcal{H}(\Phi,\Delta u,\Lambda,t) & = \Delta u^2(t) \nonumber  + \lambda_2 \big[ -R(\theta,p)  (\Delta u(t) + u_{\rm nom}(t + \Delta t)  - p)   \big] \\
&+ \lambda_1 \big[  \omega(p) + z(\theta,p)  (\Delta u(t) + u_{\rm nom}(t + \Delta t)  - p)   -   D(\theta,p)   R(\theta,p)   (\Delta u(t) + u_{\rm nom}(t + \Delta t)  - p)  \big],
\end{align}
where $\Phi \equiv \begin{bmatrix} \theta & p \end{bmatrix}^T$ contains the state variables and $\Lambda \equiv \begin{bmatrix} \lambda_1 & \lambda_2 \end{bmatrix}^T$ represents Lagrange multipliers that force the dynamics to satisfy the dynamics of the adaptive reduction.  According to Pontryagin's  minimum principle \cite{kirk98}, the control that minimizes the cost functional $C$ will minimize the Hamiltonian for all admissible $\Delta u(t)$ with the dynamics of the state variables and Lagrange multipliers subject to
\begin{align}
\dot{x} &=  \frac{\partial \mathcal{H}}{\partial \Lambda},  \label{xdoteqs} \\
\dot{\Lambda} &= - \frac{\partial \mathcal{H}}{\partial x} \label{lambdadoteqs}.
\end{align}
Evaluation of \eqref{xdoteqs} returns the state equations of the adaptive reduction from \eqref{adaptiveframe}.  Evaluation of \eqref{lambdadoteqs} yields
\begin{align}  \label{adaptivelambdas}
\dot{\lambda}_1 &=  -\lambda_1 \big[ z_\theta (\Delta u(t) + u_{\rm nom}(t + \Delta t) - p) - (D_\theta R + D R_\theta) (\Delta u(t) + u_{\rm nom}(t + \Delta t) - p)   \big] \nonumber \\
& \quad + \lambda_2 R_\theta  (\Delta u(t) + u_{\rm nom}(t + \Delta t) - p), \nonumber \\
\dot{\lambda}_2 &= -\lambda_1 \big[  \omega_p + z_p  (\Delta u(t) + u_{\rm nom}(t + \Delta t) - p) - z - (D_p R + D R_p) (\Delta u(t) + u_{\rm nom}(t + \Delta t) - p) + DR  \big] \nonumber \\
& \quad - \lambda_2 \big[R-  R_p  (\Delta u(t) + u_{\rm nom}(t + \Delta t) - p)  \big],
\end{align}
where, for instance, the notation $z_\theta$ corresponds to the partial of $z$ with respect to $\theta$ evaluated at both $\theta$ and $p$.  Noting that the Hamiltonian \eqref{hameq} is quadratic in $\Delta u$, the admissible control $\Delta u(t)$ that minimizes the Hamiltonian at any given moment is 
\begin{equation} \label{adaptiveu}
\Delta u = \min \left(  \max \left( \frac{-z(\theta,p) \lambda_1 + \lambda_1 D(\theta,p) R(\theta,p) + R(\theta,p) \lambda_2 }{2}, \Delta u_{\rm min} \right), \Delta u_{\rm max} \right).
\end{equation}
In the adaptive reduction framework, letting $\theta_{\rm ent}(t)$ and $p_{\rm ent}(t)$ correspond to the $\theta$ and $p$ values on the entrained periodic orbit, the boundary conditions prescribed by the control problem are $\theta(t_0) = \theta_{\rm ent}(t_0)$, $\theta(t_0 + T_f) = \theta_{\rm ent}(t_0 + \Delta t)$, $p(t_0) = p_{\rm ent}(t_0)$, and  $p(t_0 + T_f) = p_{\rm ent}(t_0 + \Delta t)$.  In order to solve this two-point boundary value problem, it is necessary to identify the correct choice of $\lambda_1(t_0)$ and $\lambda_2(t_0)$ that yield the prescribed final conditions.  This is accomplished in this work by providing an initial guess and subsequently updating the initial values of these Lagrange multipliers using a Newton iteration until convergence is achieved.   

\subsection{Two Alternative Optimal Control Formulations}
In addition to the optimal control formulation that uses the adaptive phase-amplitude reduction, two additional control formulations used for comparison purposes.  The first control strategy uses the first order accurate phase-isostable reduced equations from \eqref{isored} with dynamics that follow
\begin{align} \label{nonadaptive}
\dot{\theta} &= \omega + z(\theta) (\Delta u(t) + u_{\rm nom}(t_s)), \nonumber \\
\dot{\psi} &= \kappa \psi + i(\theta)  (\Delta u(t) + u_{\rm nom}(t_s)),
\end{align}
where $\omega$ and $\kappa$ are evaluated at $p_0$, $U(t) = \delta u(t)$ as defined directly above \eqref{adaptivemod}, with  $z(\theta) = Z(\theta,p_0)  \cdot \delta$ and $i(\theta) = I(\theta,p_0) \cdot \delta$ denoting the effective phase and isostable response curves.  Compared with the adaptive reduction \eqref{adaptiveframe}, Equation \eqref{nonadaptive} does not allow for adjustments of the underlying system parameter $p$, but rather, explicitly takes into account the isostable coordinates defined in reference to the periodic orbit $x^\gamma_{p_0}$.  Such a formulation allows for the definition of a cost functional, similar to the one proposed in \cite{mong19b}, that balances the tradeoff between the magnitude of the isostable coordinate and the control input
\begin{equation} \label{cost2}
C_2 = \int_{t_0}^{t_0 + T_f} \Delta u^2(t) + \beta \psi^2(t) dt,
\end{equation}
where $\beta$ is a positive constant.  Using the same control objective, same limits on the optimal control, and same initial and final conditions as those from the previous sections, one can follow an identical set of arguments to write the Hamiltonian associated with the cost functional \eqref{cost2} 
\begin{align} \label{ham2}
\mathcal{H}_2(\Phi,\Delta u,\Lambda,t) & = \Delta u^2(t)  + \beta \psi^2 +   \lambda_1 \big[   \omega + z(\theta) (\Delta u(t) + u_{\rm nom}(t + \Delta t))  \big] \nonumber \\
& \quad+ \lambda_2 \big[  \kappa \psi + i(\theta)  (\Delta u(t) + u_{\rm nom}(t + \Delta t))  \big] ,
\end{align}
where $\Phi = \begin{bmatrix} \theta & \psi \end{bmatrix}^T$. Associated optimal trajectories must satisfy \eqref{nonadaptive} along with
\begin{align}  \label{phaseamplitudecontrol}
\dot{\lambda}_1 &= -\lambda_1 z_\theta(\Delta u(t) + u_{\rm nom}(t + \Delta t))  - \lambda_2 i_\theta (\Delta u(t) + u_{\rm nom}(t + \Delta t)) , \nonumber \\
\dot{\lambda}_2 &= -\lambda_2 \kappa - 2 \beta \psi, \nonumber \\
\Delta u &=  \min \left(  \max \left(  \frac{-\lambda_1 z(\theta) - \lambda_2 i(\theta) }{2}  ,  \Delta u_{\rm min}  \right), \Delta u_{\rm max} \right).
\end{align}
In the phase-isostable reduced framework, letting $\theta_{\rm ent}(t)$ and $\psi_{\rm ent}(t)$ correspond to the $\theta$ and $p$ values on the entrained periodic orbit, the boundary conditions prescribed by the control problem are $\theta(t_0) = \theta_{\rm ent}(t_0)$, $\theta(t_0 + T_f) = \theta_{\rm ent}(t_0 + \Delta t)$, $\psi(t_0) = \psi_{\rm ent}(t_0)$, and  $\psi(t_0 + T_f) = \psi_{\rm ent}(t_0 + \Delta t)$.

The second alternative control formulation considered only uses the phase dynamics, completely ignoring all amplitude coordinates.  Such a control strategy only considers the top equation from \eqref{nonadaptive} with a cost function 
\begin{equation} \label{cost3}
C_3 = \int_{t_0}^{t_0 + T_f} \Delta u^2(t) dt.
\end{equation}
The corresponding Hamiltonian is 
\begin{equation}
\mathcal{H}_3(\theta,\Delta u,\lambda_1,t) = \Delta u^2(t) +  \lambda_1 \big[   \omega + z(\theta) (\Delta u(t) + u_{\rm nom}(t + \Delta t))  \big],
\end{equation}
with associated optimal control trajectories that must satisfy
\begin{align} \label{phaseonlyred}
\dot{\theta} &= \omega + z(\theta) (\Delta u(t) + u_{\rm nom}(t + \Delta t)), \nonumber \\
\dot{\lambda}_1 &= -\lambda_1 z_\theta(\theta) (\Delta u(t) + u_{\rm nom}(t + \Delta t)), \nonumber \\
\Delta u &= \min \left(  \max \left(-  \lambda_1 z(\theta)/2  ,  \Delta u_{\rm min}  \right), \Delta u_{\rm max} \right).
\end{align}
Without  information about the isostable dynamics, the boundary conditions only involve the phase coordinates and are taken to be $\theta(t_0) = \theta_{\rm ent}(t_0)$, $\theta(t_0 + T_f) = \theta_{\rm ent}(t_0 + \Delta t)$.  Related optimal control problems that only account for the phase dynamics were considered in \cite{moeh06} and \cite{nabi132}.

\section{Optimal Control Application in a Simple Model for Circadian Oscillations} \label{simpleresult}
Results are presented below that use adaptive phase-isostable reduction framework in conjunction with the prototype optimal control problem.  In the context of circadian cycles, the optimal control problem can be viewed as a strategy that seeks to mitigate circadian misalignment (i.e.,~jet-lag) caused by a sudden shift in the environmental time occurring as a result of rapid travel through multiple time zones.  In this section, a simple three-dimensional model for circadian oscillations is considered.  In the section to follow, a collective oscillation that represents the aggregate behavior of 3000 coupled oscillators will be considered. 


For this example, consider a simple three-dimensional model for circadian oscillations \cite{gonz05}
\begin{align} \label{circmodel}
\dot{B} &= v_1 \frac{K_1^n}{K_1^n + D^n} - v_2 \frac{B}{K_2+B} +  L_c + L_{\rm nom}(t_s) + \Delta L(t),\nonumber \\
\dot{C} &= k_3 B - v_4 \frac{C}{K_4+C}, \nonumber \\
\dot{D} &= k_5 C - v_6 \frac{D}{K_6+D}.
\end{align}
In the above equations, $B$, $C$, and $D$, represent concentrations of mRNA of the clock gene, associated protein, and nuclear form of the protein, respectively.    The term $L_{\rm nom}(t_s)$ represents a nominal 24-hour light-dark cycle taken to be
\begin{equation} \label{lightvalue}
L_{\rm nom}(t_s) = L_0  \left[  \frac{1}{1+\exp(-5 (t_{\rm s} -6))}  - \frac{1}{1+\exp(-5 (t_{\rm s} -18))}   \right],  \\
\end{equation}
where $L_0$ is the maximum uncontrolled light intensity and $t_s = {\rm mod}(t + \Delta t,24)$ with $\Delta t$ being a time shift.   The term $\Delta L(t)$ is taken to be a control input with the overall light intensity bounded by $L_{\rm min} \leq L_{\rm nom}(t_s) + \Delta L(t) \leq L_{\rm max}$.  Matching the formulation given in \eqref{adaptivemod},  $U(t) = \delta (\Delta L(t) + u_{\rm nom}(t))$ where $\delta = \begin{bmatrix} 1 & 0 & 0 \end{bmatrix}^T$.   Remaining parameters are taken to be $n = 6$, $v_1 = 0.84, v_2 = 0.42, v_4 = 0.35, v_6 = 0.35, K_1 = 1, K_2 = 1, K_4 = 1, K_6 = 1, k_3 = 0.7$, and  $k_5 = 0.7$.  With this parameter set, when $L_c = L_{\rm nom} = \Delta L = 0$, (i.e.,~in the absence of light) the resulting stable limit cycle has a period of 24.2 hours.  

\begin{figure}[htb]
\begin{center}
\includegraphics[height=3.0 in]{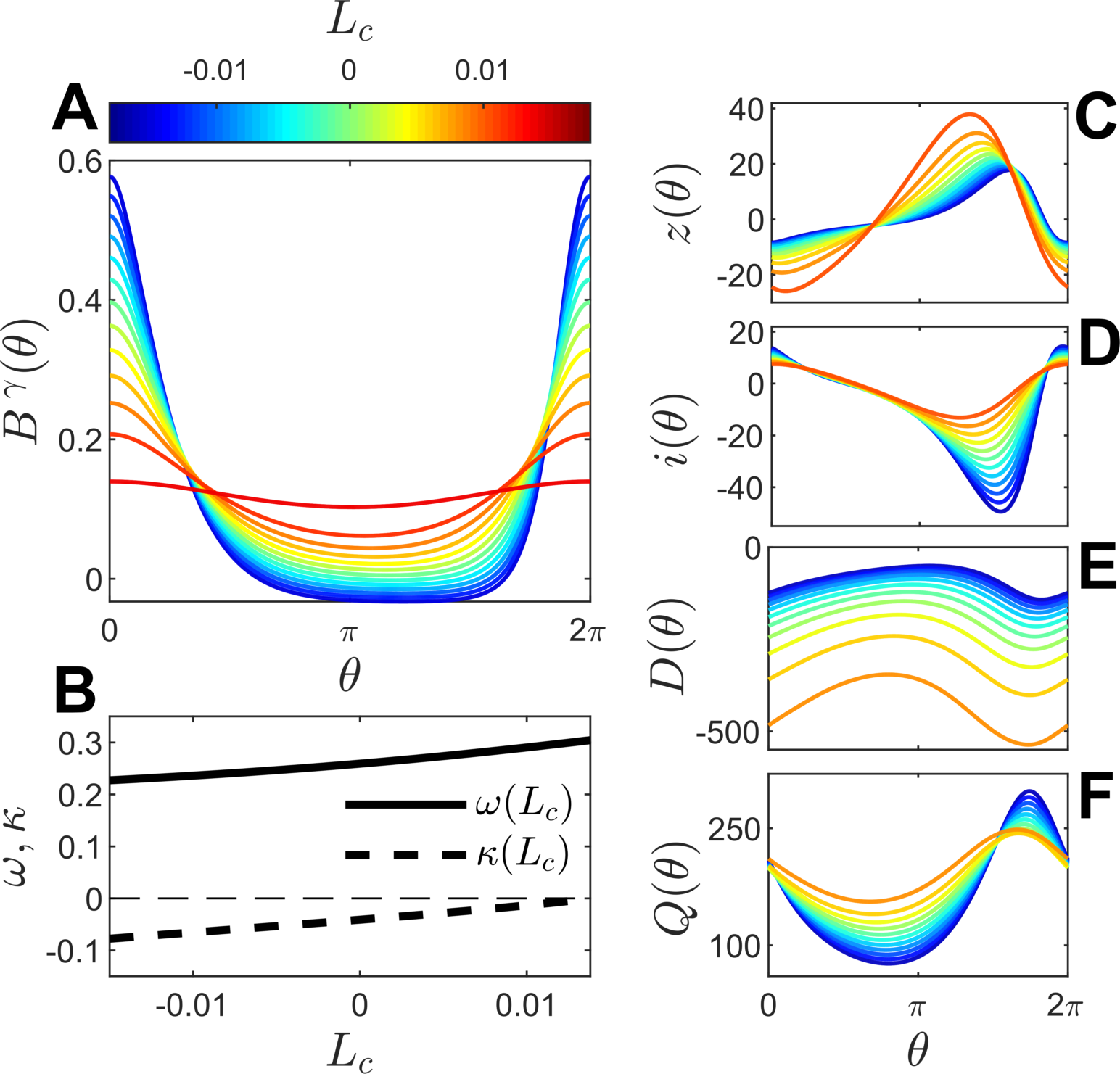}
\end{center}
\caption{  Traces of $B^\gamma(\theta)$ are shown for \eqref{circmodel} using various constant values of $L_c$.  The resulting natural frequencies and Floquet exponents are shown in panel B.  The dashed line at $\kappa = 0$ is shown for reference, emphasizing that all considered orbits are stable.  Phase and isostable response curves used in the adaptive reduction are shown in panels C and D, respectively.  The terms $D$ and $Q$ characterize how changing the nominal parameter influences the phase and isostable coordinate and are shown in panels E and F, respectively.  Note that both the magnitude of oscillations as well as the magnitude of the Floquet exponents become smaller as $L_c$ increases resulting less robust oscillations. }
\label{adaptivecurves}
\end{figure}

Letting $L_c$ be the time-varying parameter used in the adaptive reduction \eqref{adaptivemod}, the traces of $B(\theta)$ on the resulting periodic orbits $x^\gamma_{L_c}$ are shown in panel A of Figure \ref{adaptivecurves}. The periodic orbits are appropriately shifted so that $\theta = 0$ corresponds to the moment that $B(\theta)$ reaches its peak value.  For all values of $L_c$ considered, the periodic orbit has only one non-negligible Floquet exponent (the other Floquet exponent is negative and large in magnitude).  Along with this Floquet exponent, the nominal period is shown as a function of $L_c$ in panel $B$.  The resulting phase and isostable response curves for various values of $L_c$ are computed using an adjoint method described in \cite{wils19complex} and shown in panels C and D, respectively.  Finally, $D$ and $Q$ are computed according to \eqref{dandq} and shown in panels E and F, respectively.

For the moment, the optimal control problem will be considered taking $L_0 = 0$.  Additionally, it will be assumed that  $L_{\rm max} = + \infty$, and $L_{\rm min} = -\infty$.  This situation represents a situation where no entraining stimulus is applied and no bounds on the magnitude of light are present.  While this situation is not realistic, it provides insight into the limitations of each optimal control formulation.  Note that when $L_0 = 0$, the system is technically not entrained to any external stimulus.  Nevertheless, the optimal control frameworks proposed in Section \ref{optsec} can still be implemented by defining $x^\gamma_{{\rm ent}} \equiv x^\gamma_{p_0}$ with  $x^\gamma_{p_0} |_{\theta = 0} = x^\gamma_{{\rm ent}}|_{t = 0}$.

For various choices of $\Delta t$, an optimal control is computed using the adaptive reduction framework (with Equations \eqref{adaptiveframe}, \eqref{adaptivelambdas} and \eqref{adaptiveu}), the first order accurate phase-amplitude reduction strategy (with Equations \eqref{nonadaptive} and \eqref{phaseamplitudecontrol}), and the phase-only reduction strategy (with Equation \eqref{phaseonlyred}).  For the phase-amplitude and phase-only reduction strategies, $L_c$ is taken to be 0.  The constant $\beta$ (which sets the penalty for large isostable coordinates) is taken to be $10^{-4}$ for the phase-amplitude control strategy.  $T_f$ is taken to be 24 hours for all control strategies.  Optimal inputs $\Delta L(t)$ are computed using each reduction strategy, with the resulting inputs applied to the full model \eqref{circmodel}. 

  \begin{figure}[htb]
\begin{center}
\includegraphics[height=2.5 in]{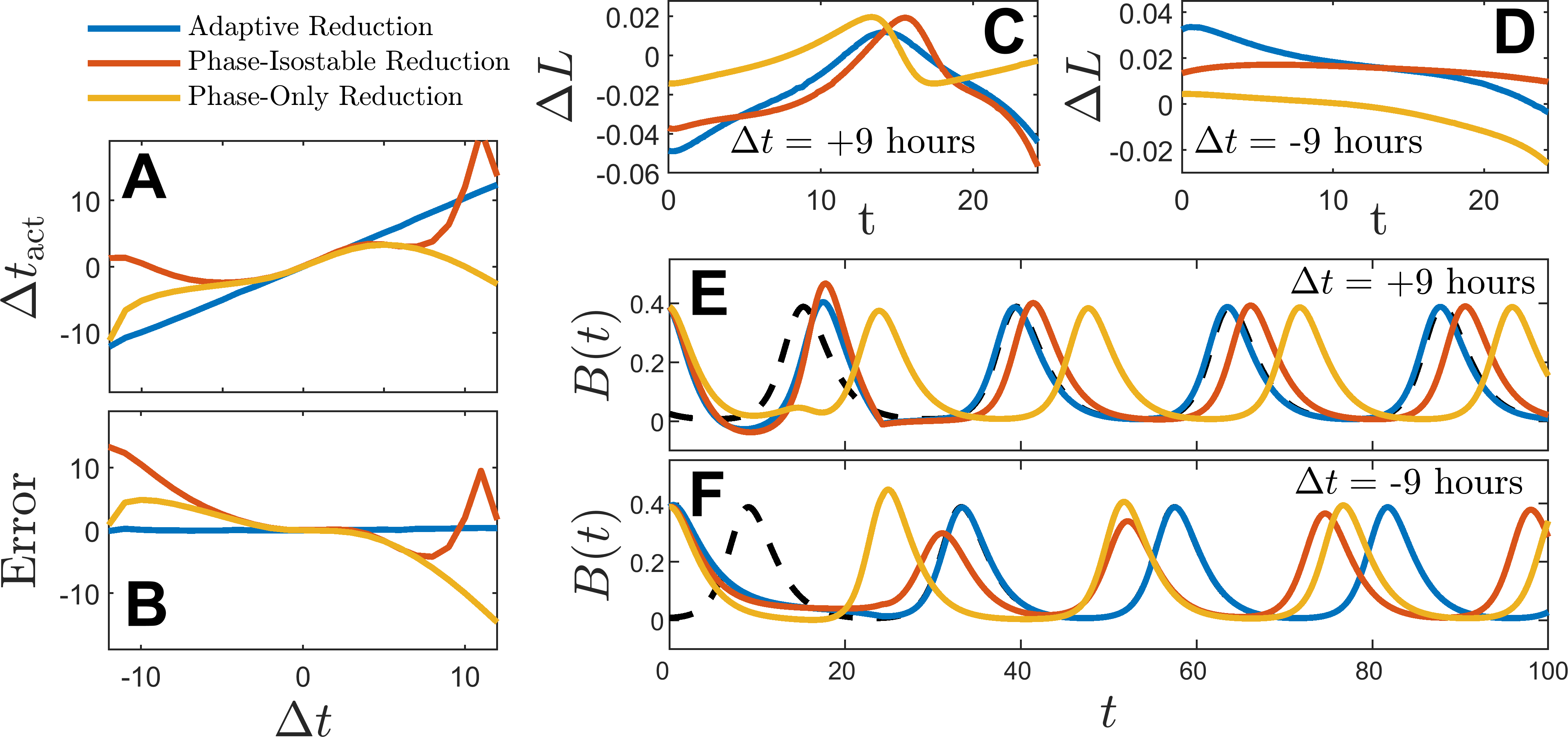}
\end{center}
\caption{ Optimal control results when no external entraining stimulus is applied and no constraints are placed on $\Delta L(t)$.   The resulting optimal controls computed using the adaptive reduction, phase-isostable reduction, and the phase-only reduction are applied to the full model equations \eqref{circmodel}.  The specified $\Delta t$ is plotted against the resulting time change, $\Delta t_{\rm act}$ in Panel A.  The error, defined to be $\Delta t_{\rm act} - \Delta t$ is shown in panel B.  Specific inputs are shown for $\Delta t = +9$ and -9 hours in panels C and D using the specified method. Resulting traces of $B(t)$ when these inputs are applied to the full model are shown in Panels E and F, respectively.   The black dashed line shows the fully entrained solution after the time shift is applied.  While all methods perform reasonably well for for $|\Delta t| < 3$ hours, the phase-isostable and phase-only optimal control strategies that yield inputs that give unexpected and inaccurate results.  Conversely, the adaptive reduction strategy accurately modifies the phase of the oscillation for any prescribed $\Delta t$.  }
\label{noconstraintresult}
\end{figure}

 Resulting time shifts are shown in panel A of Figure \ref{noconstraintresult}.  Note that because there is no external entraining stimulus in this example, mandating a time shift $\Delta t$ is equivalent to mandating a phase shift $\Delta \theta = 2 \pi \Delta t/24.2$, and $t_{\rm act}$ is the actual resulting phase shift when the input is applied to the full model \eqref{circmodel}.  $\Delta t_{\rm act}$ is computed by considering the infinite-time behavior after the input is applied.  In panel B, ${\rm Error} = \Delta t_{\rm act} - \Delta t$ is shown.  For small values of $\Delta t$, all reduction frameworks produce stimuli that achieve the control objective.  As the required $\Delta t$ shift increases, errors tend to increase for the phase-only and the phase-amplitude control strategies.  By contrast, the adaptive reduction strategy yields results that are nearly perfect, with the maximum value of $| \Delta t_{\rm act} - \Delta t|$ being less than 0.4 hours for any choice of $\Delta t$. Panel C (resp.,~D) show resulting control inputs for each strategy that result when $\Delta t = +9$ hours (resp.~-9 hours).  Panels panel E and F show corresponding traces of $B(t)$ that result when those inputs are applied to the full model \eqref{circmodel}.

Next, a situation where $L_0 = 0.01$ is considered next so that the limit cycle is nominally entrained to a the external 24-hour light-dark cycle \eqref{lightvalue}.  Noting that the nominal phase response curves from panel C of Figure \ref{adaptivecurves} can take values that are on the order of approximately 20, this input is quite large for this system. Such situations have traditionally been difficult to approach with standard phase reduction techniques (as discussed in \cite{diek16}) because the entraining input is often strong enough to drive the state far from its unperturbed limit cycle, thereby invalidating the assumptions of the phase reduction.  Alternative techniques that consider the dynamics of Poincar\'e maps \cite{diek16}, \cite{diek18} or those that consider the entrained orbit itself to be the nominal periodic orbit \cite{wils20ddred} have been useful in some situations.  As shown in the results to follow, the control strategy that uses adaptive phase-amplitude reduction yields accurate results while the phase-only and first order accurate phase-amplitude reduction strategies fail.

\begin{figure}[htb]
\begin{center}
\includegraphics[height=2.5 in]{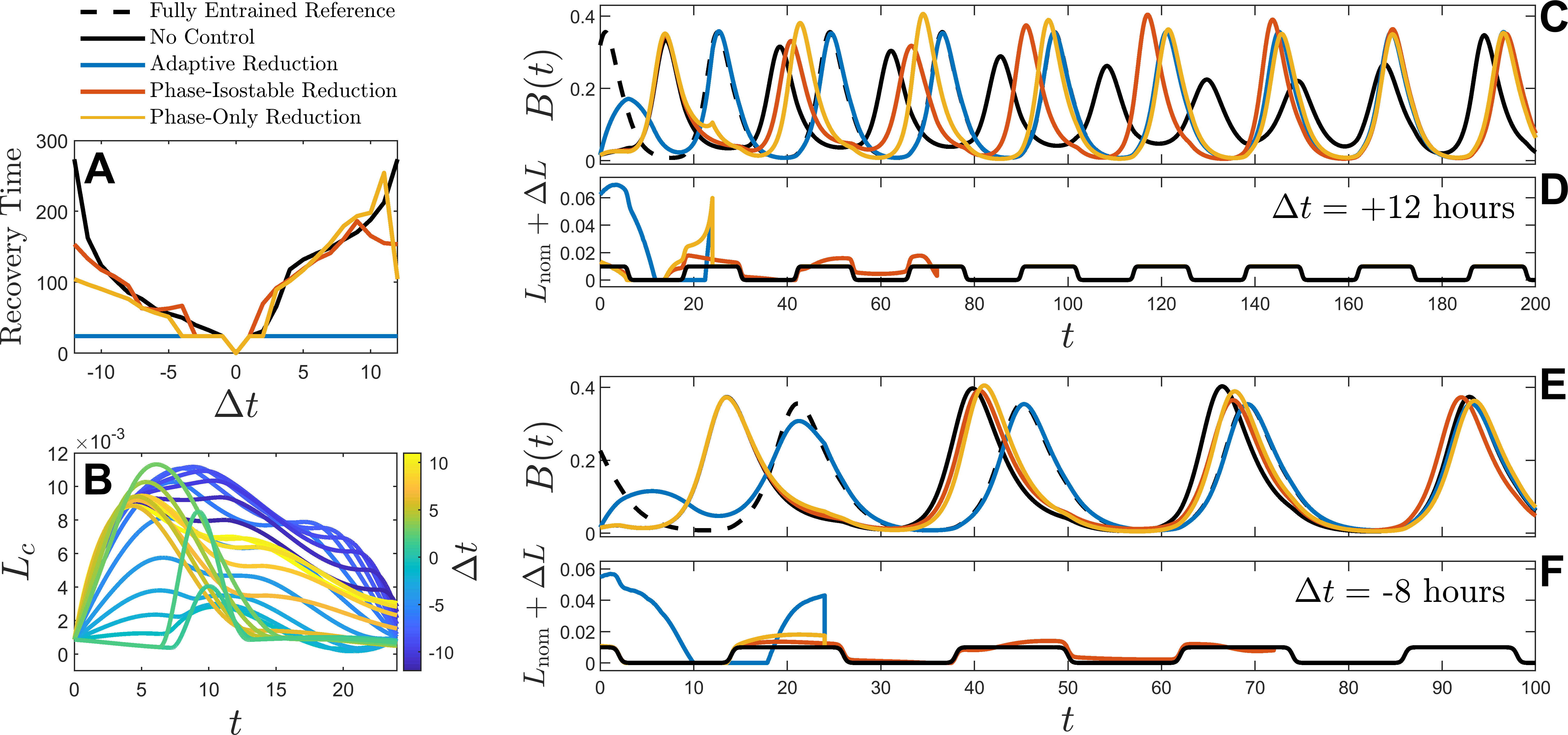}
\end{center}
\caption{Optimal Control results with a nominal external entraining stimulus and constraints on the applied stimulus are considered.   Panel A shows recovery times when applying the resulting  optimal inputs determined from each of the three optimal control formulations to the full model equations \eqref{circmodel}.  The uncontrolled recovery times are also showed for reference.  The formulation with the adaptive reduction significantly outperforms the other two control formulations.  Indeed, the phase-only and first order accurate phase-isostable formulations often yield control inputs that {\it increase} the recovery time.  Panel B shows the value of $p$, i.e.,~the parameter that corresponds to the nominal periodic orbit during the application of the optimal stimulus for various time shifts $\Delta t$.  Notice that $p$ increases to much larger values for all but the smallest magnitude choices of $\Delta t$.  These changes in $p$ represent large deviations from the nominal limit cycle which has $p_0 = 0$.  The adaptive reduction has the ability to explicitly incorporate these shifts, however,  the phase-amplitude and phase-isostable reductions cannot explicitly account for these shifts.  This fundamental difference between the reduction frameworks leads to the difference in efficacy between the resulting optimal inputs.   Panels C and D show the recovery for the optimal controls identified by each of the optimal control formulations for $\Delta t = +12$ hours. The black curve also shows the nominal recovery when $\Delta L(t) = 0$. The dashed line in panel C represents the limit cycle that is fully entrained to the shifted nominal input $L_{\rm nom}(t + \Delta t)$.  Panels E and F show the same information for a time shift of -8 hours.  }
\label{constrainedresult}
\end{figure}

For this example, $L_{\rm min} = 0$ is chosen in optimal control problem along with $L_{\rm max} = +\infty$ to mandate $L_{\rm nom} + \Delta L \geq 0$ and reflect the constraint that negative light cannot  be applied.  Once again, for various choices of $\Delta t$, an optimal control is computed using the adaptive reduction framework (from Equations \eqref{adaptiveframe}, \eqref{adaptivelambdas} and \eqref{adaptiveu}), the first order accurate phase-amplitude reduction strategy (from Equations \eqref{nonadaptive} and \eqref{phaseamplitudecontrol}), and the phase-only reduction strategy (from Equation \eqref{phaseonlyred}).  For the phase-only reduction strategy $\beta$ is set to zero, since the nominal 24-hour light-dark cycle makes it difficult to limit the magnitude of the isostable coordinate.  In all optimal control computations, $t_0 = 0$ which corresponds to the middle of subjective night.  $T_f$ is taken to be 24 hours for the phase-only and adaptive reduction strategies.  $T_f$ is taken to be 72 hours for for the phase-amplitude reduction strategy, it is not possible to find numerical solutions to the optimal control problem for smaller values of $T_f$.   Resulting inputs are then applied to the full model equations \eqref{circmodel} with results shown in Figure \ref{constrainedresult}.  Because the model is entrained to $L_{\rm nom}(t)$ in the absence of additional input, the resulting recovery time, defined to be the time it takes for the phase to return to and stay within one hour if its steady state behavior, is shown in Panel A.  Results are also compared to the nominal, uncontrolled recovery time (i.e.~taking $\Delta L = 0$).    Note that the phase is only measured at $\theta = 0$, which can be observed once per cycle when $B(t)$ reaches a local maximum.  The recovery time is inferred by computing the difference in time between the moment $\theta$ reaches 0 for the fully entrained reference and the reentraining simulation and interpolating the time difference for all values in between.  If the time difference is less than one hour the first moment $\theta = 0$ is reached after input is applied for the reentraining system, the recovery time is taken to be equal to $T_f$.

The optimal controls computed according to the phase-isostable and the phase-only reduction strategies generally have only a small influence on the recovery times and even occasionally result in worse outcomes than if no input was applied at all.  Conversely, the when the control identified when using the adaptive phase-amplitude reduction is applied to the full model equations, the recovery time exactly 24 hours indicating the system is recovered by the time the first measurement of $\theta = 0$ is made after the optimal stimulus is applied.  Panel B of Figure \ref{constrainedresult} shows $L_c$ when the optimal control is applied to the adaptive reduction.  Panels C and D (resp.,~E and F) show the optimal control inputs and traces of $B(t)$ for a time shift of $\Delta t = +12$ hours (resp.~-8 hours).    Black lines also show traces of $B(t)$ and the nominal light-dark cycle input during the uncontrolled recovery.


\section{Control of Population-Level Oscillations and Data-Driven Inference of Reduced Equations} \label{popsec}
Here, a model for the behavior of a population coupled circadian oscillators will be considered:
\begin{align} \label{popmodel}
\dot{B_i} &= v_1 \frac{K_1^n}{K_1^n + D_i^n} - v_2 \frac{B_i}{K_2+B_i}  +     h_c \frac{K F}{K_c + K F}  +  \alpha_i ( L_c   + L_{\rm nom}(t_s) + \Delta L(t) ),\nonumber \\
\dot{C_i} &= k_3 B_i - v_4 \frac{C_i}{K_4+C_i}, \nonumber \\
\dot{D_i} &= k_5 C_i - v_6 \frac{D_i}{K_6+D_i}, \nonumber \\
\dot{E_i} & = k_7 B_i - v_8 \frac{E_i}{K_8 + E_i}, \quad i = 1,\dots,N.
\end{align}
The above model considers a population of $N = 3000$ coupled oscillators of the same form as \eqref{circmodel} with the addition of the dynamics of a neurotransmitter $E_i$ that allows communication between cells.  Assuming that the spatial transmission of the neurotransmitter is fast relative to the 24-hour scale of oscillation, effective mean-field coupling is assumed with $F \equiv (1/N)\sum_{k = 1}^N E_i$, which enters into the equation for the variable $B$.  Additionally, in \eqref{popmodel} each neuron has an intrinsic sensitivity to light, $\alpha_i$, drawn from the distribution $ \alpha_i = \max(1 + 0.4 N(0,1),0)$ where $N(0,1)$ is a normal distribution with zero mean and unit variance.  Nominal parameters in this model are taken to be $n = 5, v1 = 0.868, v2 = 0.634, k_3 = 0.7, v_4 = 0.35, k_5 = 0.7, v_6 = 0.35,k_7 = 0.35, v_8 = 1, K_1 = 1, K_2 = 1  , K_4 = 1, K_6 = 1, K_8 = 1, h_c = 1, K_c = 1$, and $K = 0.5$.  In order to incorporate heterogeneity in the model, the parameters $k_3, k_5, k_7, v_4,v_6$, and $v_8$ are drawn from a normal distribution with the mean being the nominal parameter value and a standard deviation equal to 0.01.  In all simulations, independent and identically distributed zero mean white noise with intensity 0.0002 is added to the variable $B_i$ for each oscillator.  Like in the single oscillator equations \eqref{circmodel},  the external light-dark cycle $L_{\rm nom}(t_s)$ in \eqref{popmodel} takes the form \eqref{lightvalue}.

 While each of the necessary terms of the reduction in \eqref{adaptivemod} could be computed numerically, a strategy for inferring the required terms from output data will be considered here, as would be the case for an experimental system. The control strategy from Section \ref{optsec} will then be applied to the resulting adaptive reduction.

\subsection{Data-Driven Methods for Computation of Terms of the Adaptive Reduction} \label{procsec}

When the underlying dynamical equations are known, it is relatively straightforward to compute the necessary terms of the adaptive reduction \eqref{adaptivemod} using the `adjoint method' \cite{brow03} and related equations from \cite{wils20highacc} as described in Section \ref{backsec}.  However, when the the dynamical equations are unknown, the terms of the adaptive reduction must inferred from data.   The `direct method'  \cite{izhi07}, \cite{neto12} is a well-established technique for identifying the phase response curve from \eqref{standardphase}.   Previous work \cite{wils20ddred} developed a data-driven technique for computing the terms of the phase-amplitude reduction from \eqref{isored}.  This strategy will be expanded here to compute the necessary terms of the adaptive reduction from \eqref{adaptivemod}.  Here, it is assumed that the input also the sole time-varying parameter.  While \eqref{adaptivemod} has only one isostable coordinate, this technique is relatively straightforward to extend to situations where multiple isostable coordinates are involved.

For concreteness, the model \eqref{popmodel} will be used assuming light perturbations $\Delta L(t)$ can be given and that the single output $F(t)$ (i.e.,~the mean value of $E_i$ for the population) can be measured.  For the given parameter set it will be assumed that Equation \eqref{popmodel} has only one dominant isostable coordinate.  Additionally, $L_c$ will be taken as the time-varying paramater.  In this situation, it is possible to write the equations for the adaptive phase-amplitude reduced model in the form \eqref{adaptivemod}.  The terms of the adaptive reduction will be inferred using the procedure detailed below.

\begin{figure}[htb]
\begin{center}
\includegraphics[height=2.5 in]{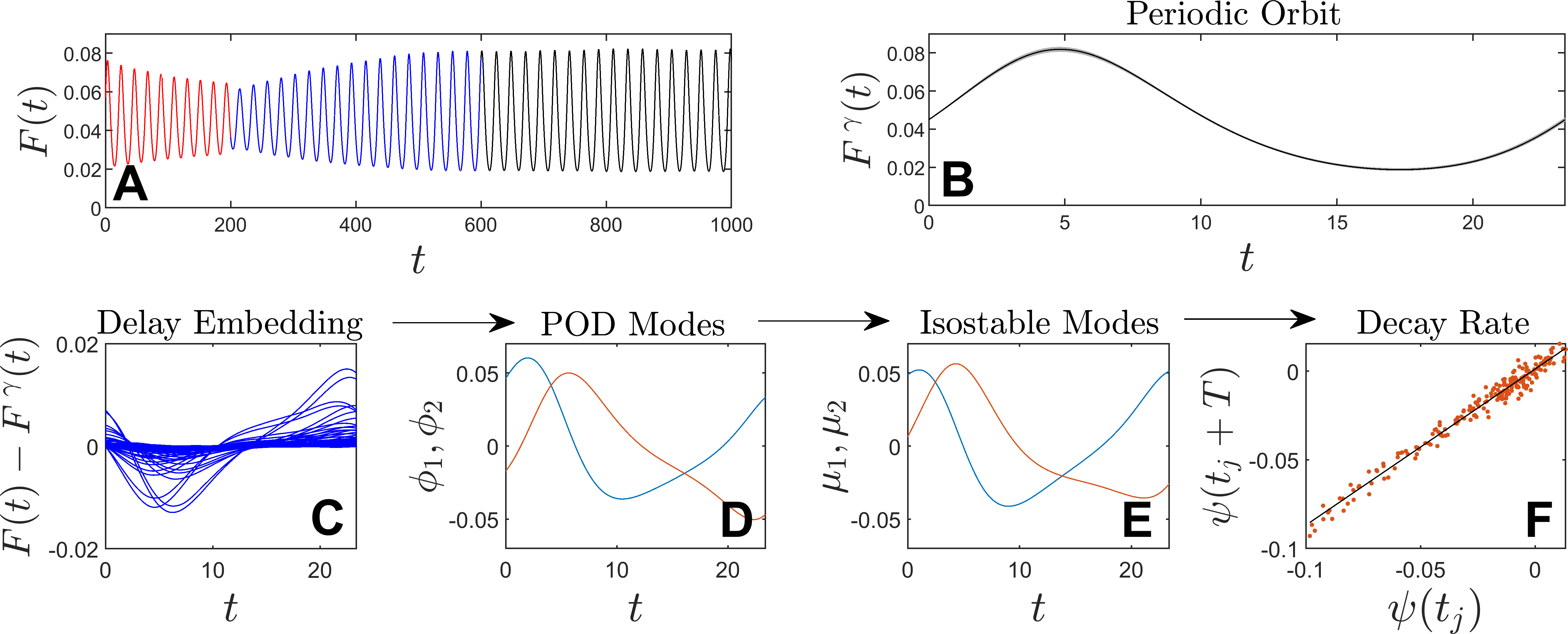}
\end{center}
\caption{A data-driven strategy for identifying phase and isostable coordinates from the large circadian model \eqref{popmodel}.  In panel A, the model is perturbed from its limit cycle by decreasing $K$ by a factor of 2 for 200 hours (red curve).  The subsequent recovery takes approximately 400 hours (blue curve) and the black curve represents the behavior once the model has recovered to the stable limit cycle.  As shown in panel B, the periodic orbit $F^\gamma(t)$ (black line) is taken to be the average of multiple cycles once transient behavior has died out.  Using data obtained from the blue portion of panel A, multiple delayed embeddings are obtained and shown in panel C that represent the deviation from the periodic orbit on a particular cycle.  POD is performed on the data from panel C and identifies two characteristic modes.  The POD modes are transformed to represent the deviation from the periodic orbit in a basis of phase and isostable coordinates with modes shown in panel E.  As described in the text, this procedure allows one to measure the isostable coordinates at specific times, and a the slope of the linear regression shown in Panel F represents the Floquet multiplier corresponding to the isostable coordinate.}
\label{datadrivenisostables}
\end{figure}

 \vskip .1 in
\noindent {\bf Step 1:}  With a static choice of $L_c$, and taking $L_{\rm nom}(t_s) + \Delta L(t) = 0$ for all time, the procedure introduced in \cite{wils20ddred} can be used to identify $\omega(L_c)$, $\kappa(L_c)$.  This strategy can also be used to identify $\psi(t_1)$ and $\theta(t_1)$ provided $\theta(t_1)\approx 0$ using a delayed embedding of the output from $F(t_1)$ to $F(t_1 + T)$.  This strategy is summarized in Figure \ref{datadrivenisostables}.  To apply this strategy, the model is first simulated until transient behaviors have died out.  A threshold is chosen to denote $\theta \approx 0$ (in this case, $\theta \approx 0$ corresponds to the moment $F(t)$ crosses 0.045 with a positive slope).   Over multiple oscillations, the average period at which this threshold is crossed is taken to be the period of oscillation, $T$.  A set of delay embeddings, each which start the moment that $F(t)$ crosses 0.045 and end $T$ hours later, is recorded and the average value of the output $F$ over these measurements is used as an approximation for the stable periodic orbit $F^\gamma(t)$ as shown in panel B of Figure \ref{datadrivenisostables}.  Once the periodic orbit is obtained, the recovery to the periodic orbit from a perturbed initial condition is considered in order to determine the isostable coordinates.  As shown in panel A,  the coupling strength $K$ is decreased 50 percent for $t \in [0,200]$ hours.  For $t > 200$, $K$ is set back to its nominal value, and a series of delayed embeddings from $t_j$ to $t_j + T$ are taken, where each $t_j$ is chosen so correspond to a moment that $F(t_j)$ crosses 0.045 with a positive slope.   This procedure is repeated over multiple trials.  The periodic orbit is subtracted from the resulting embeddings with individual traces shown in panel C.  To proceed, let $y_j \in \mathbb{R}^q$ correspond to the $j^{\rm th}$ delay embedding from panel C.  In other words, $y_j$ is comprised of  the output from $t_j$ to $t_j + T$ after $F^\gamma$ is subtracted.  Here, $q = (T/\delta t) + 1$ where $\delta t$ is the sampling rate.

Proper orthogonal decomposition (POD) \cite{benn15}, \cite{holm96}, \cite{tair17} is applied to the data shown in panel C in order to identify a small number of representative modes.  In this case, 2 modes are sufficient to represent the data from panel C, denoted by $\phi_1 \in \mathbb{R}^q$ and $\phi_2 \in \mathbb{R}^q$ and shown in panel D as blue and red curves, respectively.    The transformation $\begin{bmatrix}  \mu_1  & \mu_2 \end{bmatrix} = \begin{bmatrix}  \phi_1  & \phi_2 \end{bmatrix} A$ where $A$ is a $2 \times 2$ nonsingular square matrix and is chosen so that $\mu_1 = \begin{bmatrix}  \phi_1  & \phi_2 \end{bmatrix} ^\dagger g^\gamma$ where $g^\gamma \in \mathbb{R}^q$ corresponds to $d F^\gamma / dt$ taken at the sampling rate $\delta t$ and $^\dagger$ represents the Moore-Penrose pseudoinverse.  Plots of $\mu_1$ and $\mu_2$ are shown in panel E in blue and red, and are proportional to the shifts in the output resulting from a shift in $\theta$ and $\psi$, respectively.  Furthermore as described in \cite{wils20ddred}, the vectors $\eta_1$ and $\eta_2$ defined according to, $A^{-1}   \begin{bmatrix}  \phi_1  & \phi_2 \end{bmatrix}^T = \begin{bmatrix} \eta_1 & \eta_2 \end{bmatrix}^T$ can be used to compute the phase and isostable coordinates at specific instances in time according to the relationships
\begin{align} \label{cfns}
c  \eta_1 ^T y_j &= \theta(t_j), \nonumber \\
\eta_2^T y_j &= \psi(t_j),
\end{align}
where $c$ is a constant that can be determined using a strategy discussed in \cite{wils20ddred}.  This information can be used to track $\psi$ over successive periods during the recovery as plotted in panel F.  The slope of a linear regression of this data gives an approximation of the Floquet multiplier, $\Lambda$, corresponding to $\psi$ and the Floquet exponent can be computed according to the relationship $\kappa = \log(\Lambda)/T$ .

Note that the procedure from Step 1 represents an implementation of the procedure given in Section IV,A from \cite{wils20ddred} for a constant choice of $L_c$.  The interested reader is referred to \cite{wils20ddred} for a more detailed description of this method.  This strategy is repeated to determine the associated terms of the phase and isostable coordinate reduction for various values of $L_c$.  

\begin{figure}[htb]
\begin{center}
\includegraphics[height=3.0 in]{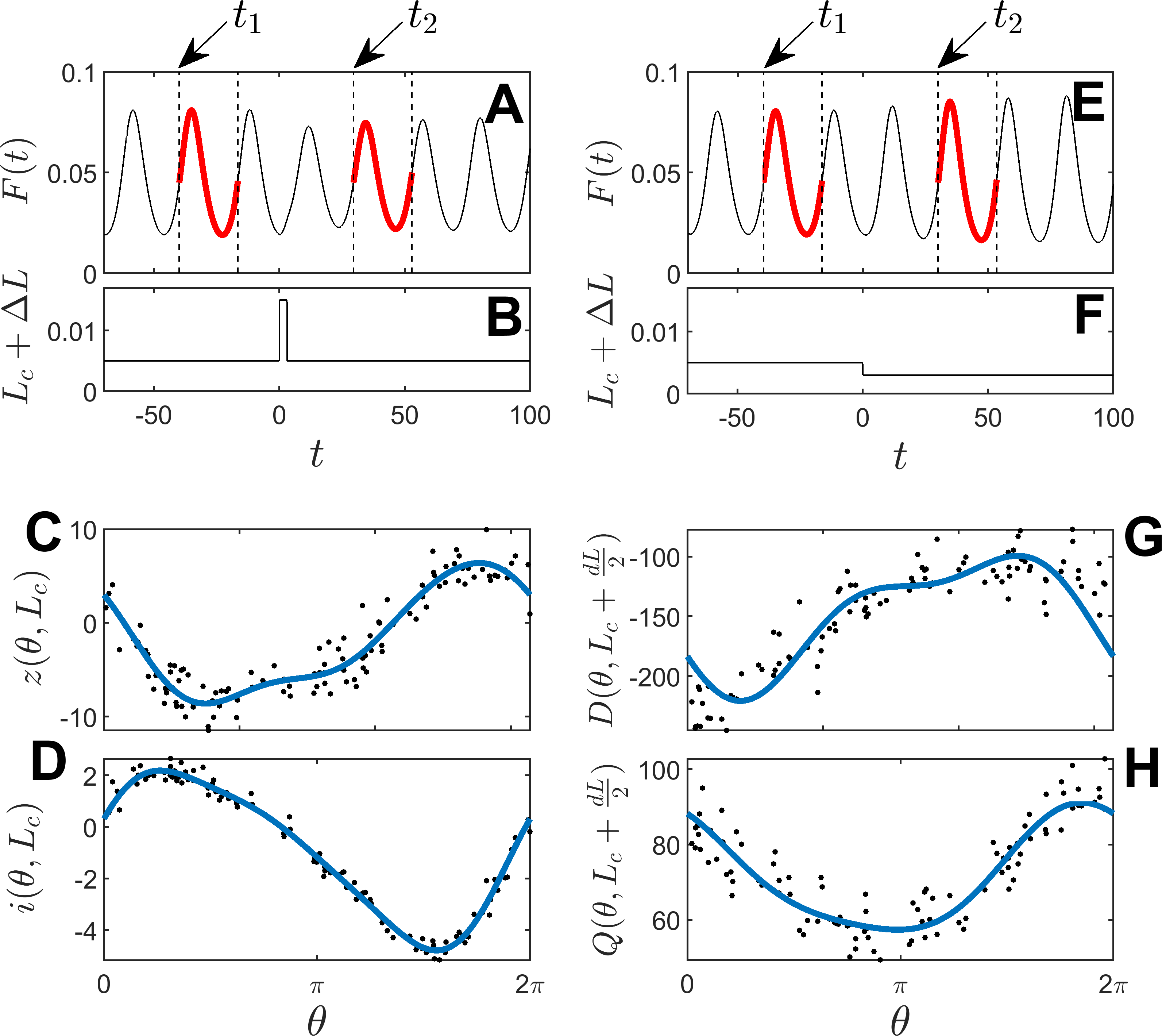}
\end{center}
\caption{Inferring the phase and isostable response to pulse and step function inputs to infer the terms of the reduction \eqref{adaptivemod}.   For a nominal value of $L_c$, a short pulse can be applied as in panel B.  Data from a complete cycle occurring both before and after the perturbation (red curves in panel A) are used to identify the phase and isostable coordinates at $t_1$ and $t_2$.  This information is then used to obtain a discrete measurement of $z(\theta,L_c)$ and $i(\theta,L_c)$,  and this process is repeated multiple times.  Each trial yielding a single measurement of both $z$ and $i$ which are represented by black dots in panels C and D.  The blue lines are obtained by fitting the sinusoidal basis of the form $\sum_{n = 0}^3 [a_n \sin(n \theta) + b_n \cos(n \theta)]$; the resulting fits are taken to be the response curves.  Likewise, panels E through H   illustrate a similar procedure for obtaining $D$ and $Q$.  Instead of applying pulses, step function inputs are applied (panel F), complete cycles before and after the step function are extracted (panel E) and used to infer the change in phase and isostable coordinates resulting from the input.  Each trial gives a single data point for $D$ and $Q$ (black dots in panels G and H).  This process is repeated for multiple trials and the same sinusoidal basis is used to fit the data (blue lines).}
\label{calculateresponseterms}
\end{figure}

 \vskip .1 in
\noindent {\bf Step 2:} For a given choice of $L_c$ the terms of $z$ and $i$ can be inferred directly by applying a series of pulse inputs and measuring the resulting change in the phase and isostable coordinates.  This can be accomplished with a strategy akin to the direct method \cite{izhi07}, \cite{neto12} whereby a pulse input of magnitude $dL$ is applied for a short duration of time as shown in panel $B$ of Figure \ref{calculateresponseterms}.  Complete cycles of the output starting at $\theta \approx 0$ (corresponding to a positive crossing $F(t) = 0.045$ threshold)  beginning at $t_1$ and $t_2$ preceding and following the input can be extracted as shown in panel A,  and the relationships from \eqref{cfns} can be used to determine $\theta(t_1)$, $\theta(t_2)$, $\psi(t_1)$, and $\psi(t_2)$.  Assuming the pulse input was applied at $t = 0$, the shift in phase can be computed by recalling that the phase would have simply increased at the rate $\omega$ had the input not been applied.  The deviation, $\Delta \theta$, from this expectation can be assumed to be the change in phase caused by the input.  Likewise, by assuming that $\psi$  decays exponentially towards zero at a rate governed by $\kappa$ in the absence of input, $\psi$ can be computed immediately before and after the applied input to infer the change in isostable coordinate $\Delta \psi$.  For the chosen parameter $L_c$, $z(\theta,L_c) \approx  \frac{\Delta \theta}{dL dt}$ and $i(\theta,L_c) \approx  \frac{ \Delta \psi}{dL dt}$, where $dt$ is the duration  of the perturbation applied at phase $\theta$.  This process can be repeated for multiple trials to obtain approximations of $i(\theta,L_c)$ and $z(\theta,L_c)$ as shown in Panels C and D of Figure \ref{calculateresponseterms} and the resulting data can be fit using an appropriate basis.

 \vskip .1 in
\noindent {\bf Step 3:}  Both $D$ and $Q$ can be inferred using a strategy similar to the one described in Step 2.  In this case, instead of a brief pulse, the total applied light is shifted from $L_c$ to $L_c + d L$ at a known time as shown in panel F of Figure \ref{calculateresponseterms}.  This shift can be viewed as a sudden change in the nominal parameter by $d L$ centered at $L_c+ dL/2$; in other words a pulse to $\dot{L}_c$ of magnitude $d L$.  Considering the structure of the reduced Equations \eqref{adaptivemod}, this pulse in $\dot{L}_c$ can be used to infer both $D$ and $Q$ in a manner similar to the direct method for identifying $z$ and $i$.   Similar to the procedure used in step 2, cycles beginning at $t_1$ and $t_2$ preceding and following the step function change in input as shown in panel E of Figure \ref{calculateresponseterms} can be isolated from the output.  Subsequently,  Equations \eqref{cfns} can be used to identify $\theta(t_1)$, $\theta(t_2)$, $\psi(t_1)$, and $\psi(t_2)$.  Note that for each trial $\theta(t_1)$ and $\psi(t_1)$ (resp.~$\theta(t_2)$ and $\psi(t_2)$) must be found using $\eta_1$ and $\eta_2$ that are obtained from Step 1 taking the nominal parameter to be $L_c$ (resp.~$L_c + dL$).  Using the same reasoning employed as part of Step 2, the change in phase, $\Delta \theta$,  and isostable coordinate, $\Delta \psi$, caused by the parameter change can be inferred yielding individual data points for $D(\theta,L_c + dL/2) = \Delta \theta / \Delta p$ and $Q(\theta,L_c + dL/2) = \Delta \psi / \Delta p$.  This process can be repeated to obtain multiple data points; curves can be fit to the resulting data as shown in panels G and H of Figure \ref{calculateresponseterms}.

\begin{figure}[htb]
\begin{center}
\includegraphics[height=2.8 in]{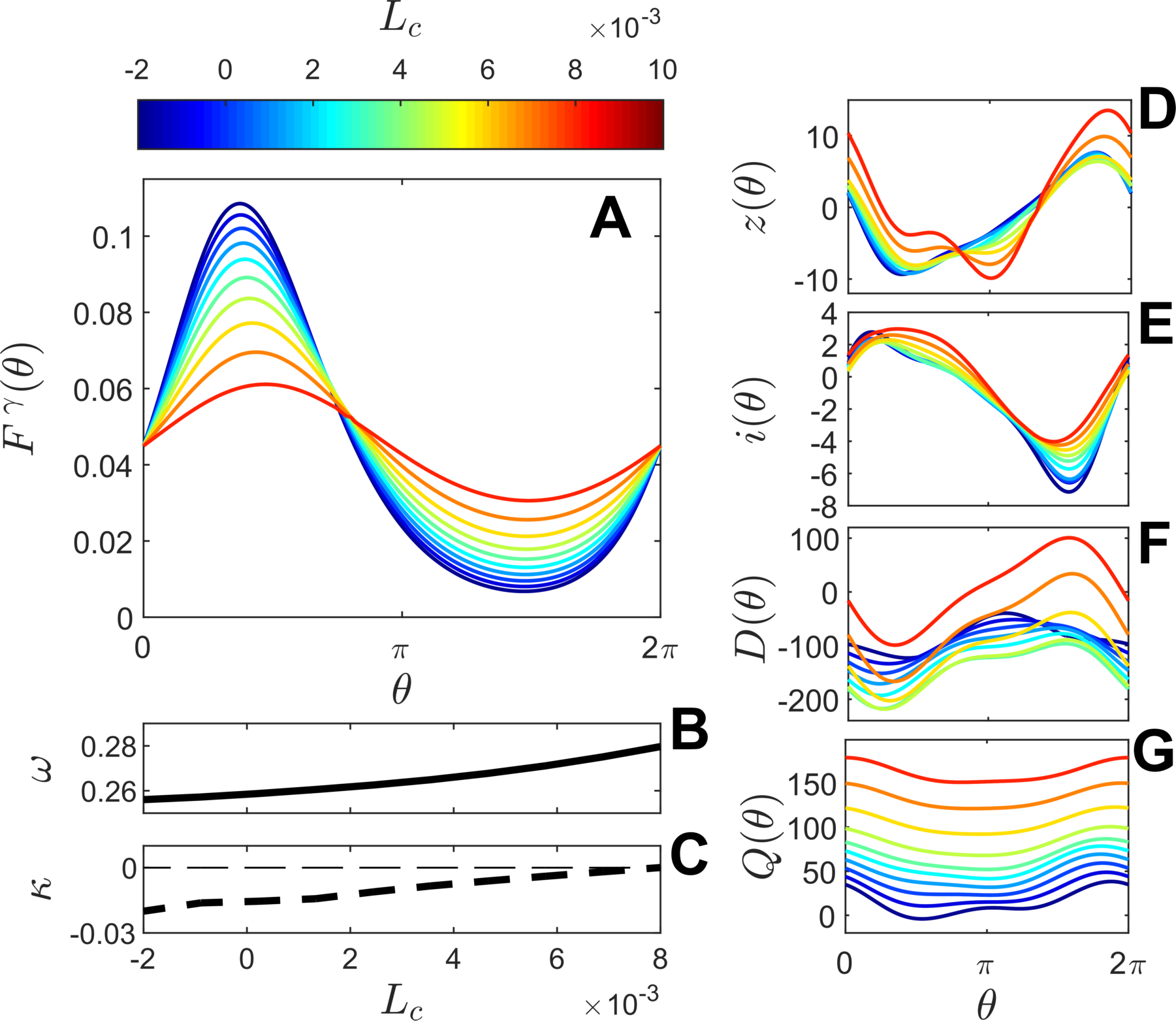}
\end{center}
\caption{ Resulting terms of the adaptive reduction \eqref{adaptivemod} resulting from the data-driven inference strategy from Section \eqref{procsec}. For all limit cycles, $\theta = 0$ corresponds to the moment that $F(t)$ crosses 0.045 with a positive slope.  As $L_c$ increases, the magnitude of the nominal oscillations decrease as shown in panel A.  Panels B and C give the natural frequency and Floquet exponent as a function of $L_c$.   Panels D-F show specific terms of the adaptive reduction for various values of $L_c$.  }
\label{reducedlarge}
\end{figure}

 \vskip .1 in
\noindent {\bf Step 4:} Steps 3 and 4 can be repeated for various choices of $L_c$ to obtain a series of curves used in the reduced equations of the form \eqref{adaptivemod}.  Resulting information obtained when applying the procedure from Steps 1 through 4 to the model from \eqref{popmodel} are shown Figure \ref{reducedlarge}.  Numerical trials are performed taking $L_c \in \{-0.01,0.01,0.03,0.05,0.07\}$, and all other curves are obtained through interpolation.

\FloatBarrier

\subsection{Optimal Control of Oscillation Timing Results}
The terms of the adaptive reduction inferred using the strategy from Section \ref{procsec} and shown in Figure \ref{reducedlarge} are used to implement the optimal control strategy for manipulating oscillation timing as described in Section \ref{optsec}.  The nominal value of $L_c$ is taken to be zero and the nominal magnitude of the light-dark cycle is taken to be $L_0 = 0.01$.  Bounds on the total allowable light input are chosen to be $L_{\rm min} = 0$ and $L_{\rm max} = \infty$ so that negative light inputs are not possible.  Recalling that the choice of $\Delta t$ in the control formulation corresponds to a sudden change in the environmental time (perhaps due to rapid travel across multiple time zones) the optimal control input is computed for various choices of $\Delta t$ with the adaptive reduction framework (from Equations \eqref{adaptiveframe}, \eqref{adaptivelambdas} and \eqref{adaptiveu}) and the first order accurate phase-amplitude reduction strategy (from Equations \eqref{nonadaptive} and \eqref{phaseamplitudecontrol}).  For when using the first order phase-amplitude reduction strategy, $\beta$ is set to zero since the nominal 24-hour light-dark cycle makes it difficult to limit the isostable coordinates.  For many choices of $\Delta t$, the phase-only optimal control equations \eqref{phaseonlyred} have solutions that grow to infinity in finite time; consequently this optimal control formulation will not be considered for this application.  In all optimal control computations, $t_0 = 0$ which corresponds to the middle of subjective night.  When possible, $T_f$ is taken to be 24 hours.  For positive time shifts and particularly large magnitude negative time shifts, solutions of the optimal control equations were not able to be found using $T_f = 24$ hours.  In these cases, $T_f$ is instead taken to be 48 hours.

\begin{figure}[htb]
\begin{center}
\includegraphics[height=2.2 in]{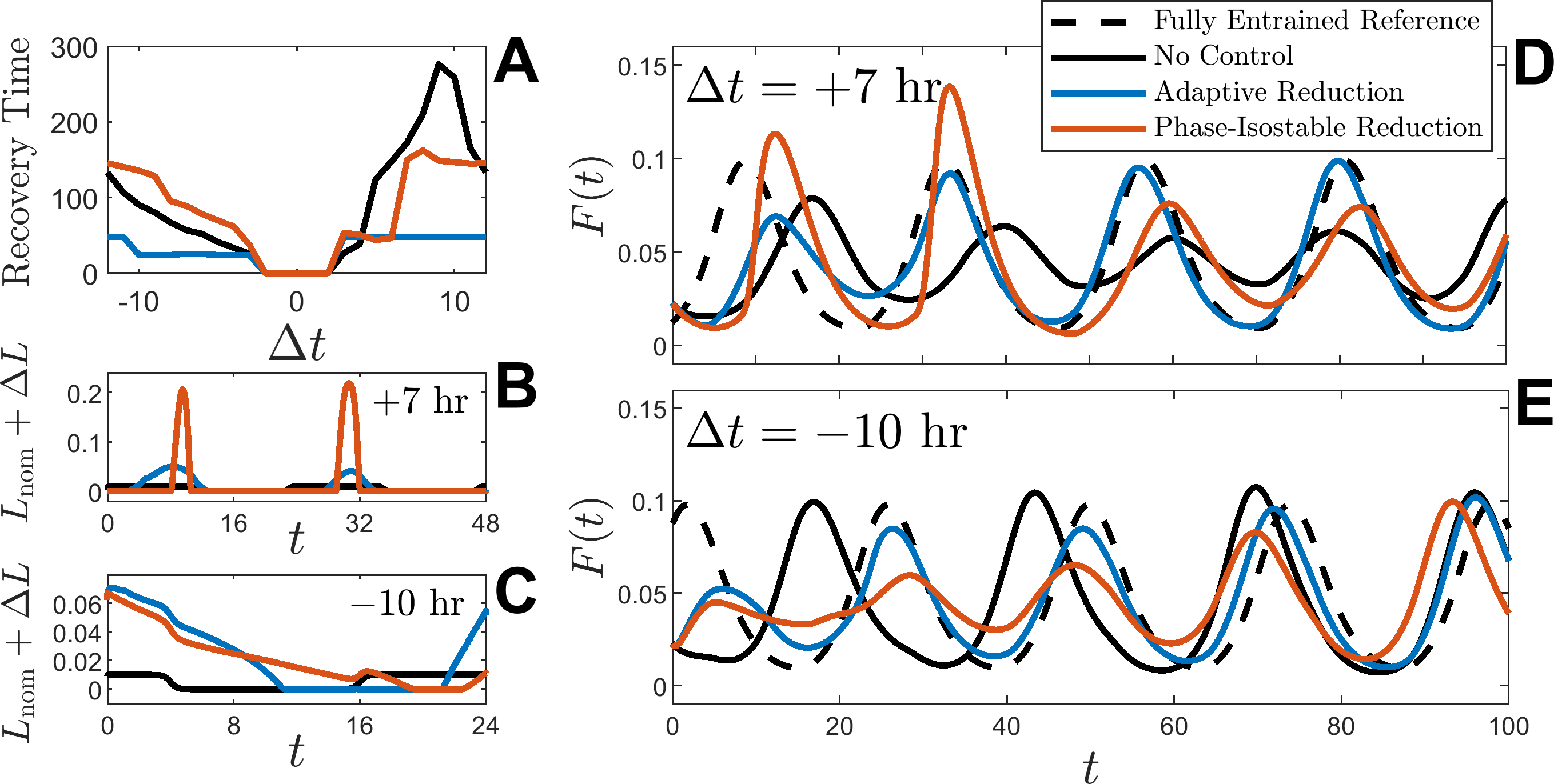}
\end{center}
\caption{Panel A shows recovery times when applying inputs obtained from the adaptive reduction and first order accurate phase-isostable reduction are shown as blue and red lines, respectively.  For reference, the black line show the uncontrolled recovery times.  Recovery is substantially hastened for optimal inputs obtained from the adaptive reduction.  Conversely, it is often worse to apply inputs obtained from the first order accurate phase-isostable reduction than it is to do nothing.  Panels B and C show example optimal control inputs for a positive and negative time shift, respectively.  The stimuli in panel B are qualitatively similar for each optimal control  both suggesting that two short pulses of light will engender the required + 7 hour phase shifts.  Nevertheless, the differences in timing and magnitude of these stimuli yield qualitatively different behaviors when applied to the full model \eqref{popmodel} as shown in panel D.  Similarly, the stimuli shown in panel C are nearly identical for  $t \in [0, 9]$, however, the differences between these stimuli lead to the markedly different recoveries shown in panel E. }
\label{recoverypopulation}
\end{figure}

The numerically determined optimal control inputs are applied to the full model equations \eqref{popmodel} with results shown in Figure \ref{recoverypopulation}.  The recovery time in panel A is defined to be the time required (after the time shift by $\Delta t$ of the light-dark cycle) for the phase to return to and stay within two hours of its steady state behavior.  The nominal recovery time that results when $\Delta L$ is zero is also shown for reference.  Note that the moment that $F(t)$ crosses 0.045 is used to correspond to the time that $\theta = 0$ and that this can only be observed once per cycle -- for this reason, if the phase is recovered at the first crossing of $\theta = 0$ after the optimal stimulus is applied, the recovery time is taken to be $T_f$.  For each optimal control determined from the adaptive reduction, recovery occurs at or very close to to $T_f$.   The optimal control derived from the first order accurate phase-isostable reduction speeds recovery time for phase advances. For phase delays, however, recovery is actually delayed when the optimal control identified from the first order phase-isostable reduction is used.    Panel B (resp.~C) shows control inputs obtained for a $\Delta t = +7$ hour (resp.~-10 hour) time shift.  The resulting model outputs when stimulus is applied to the full model equations \eqref{popmodel} are shown in panel D (resp.~E).  In both cases, the input derived from the adaptive reduction is qualitatively similar to the input from the input derived using the phase-isostable reduction, however, the quantatative differences between each input result in substantially different recovery times.   

  \begin{figure}[htb]
\begin{center}
\includegraphics[height=3.2 in]{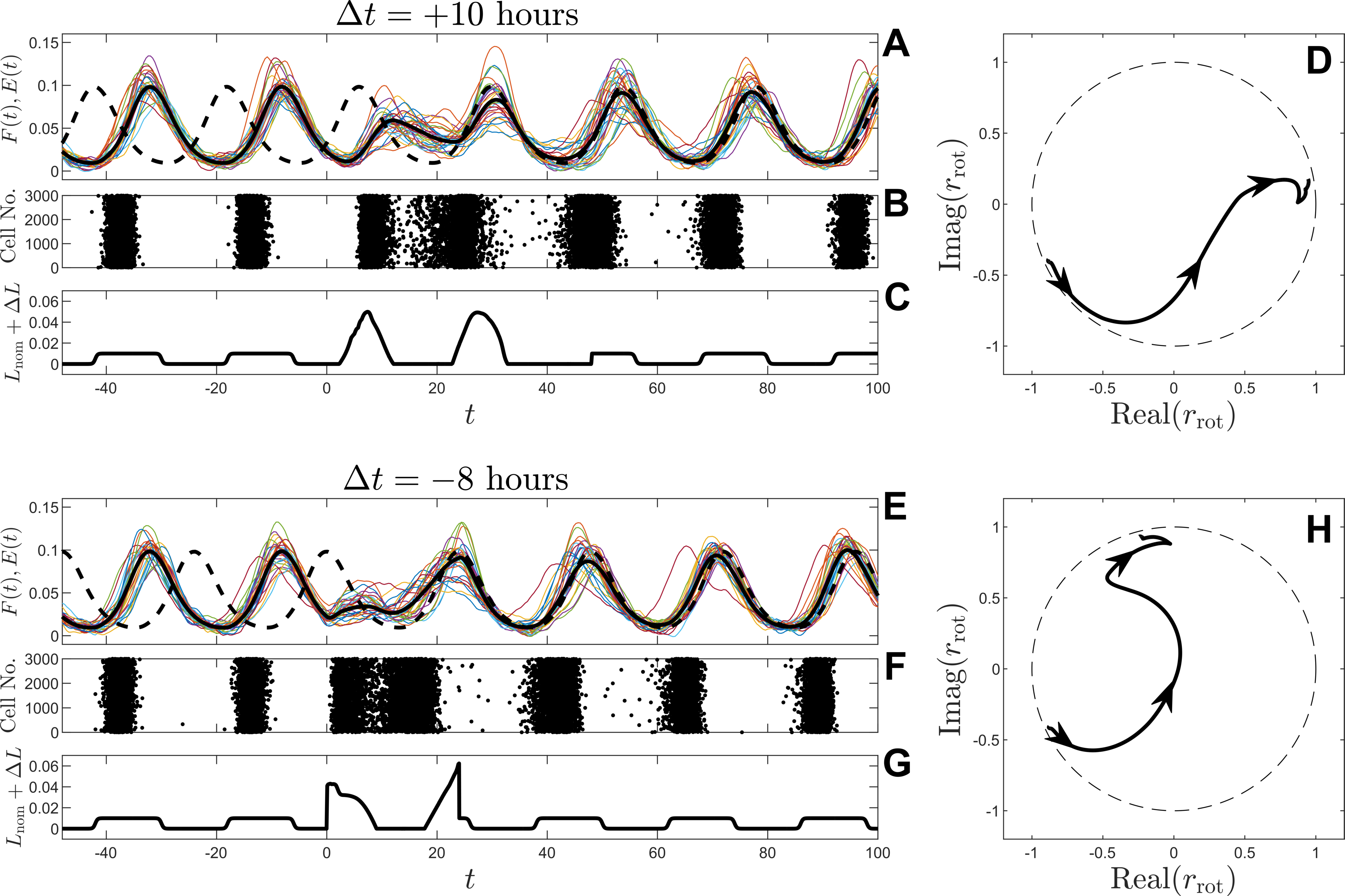}
\end{center}
\caption{ Panel A shows traces of $E_i(t)$ for 30 representative oscillators from the larger population when a time shift of $\Delta t$ = +10 hours is applied starting at $t = 0$.  The optimal control (panel C) is obtained using the adaptive reduction and is applied in order to rapidly reentrain the population oscillation to the new 24-hour light-dark cycle.  Each dot from the raster plot in panel B represents the time that $E_i$ crosses the threshold 0.035 with a positive slope taken to correspond to $\theta = 0$ for the individual oscillators.    Panel D shows the rotating order parameter from \eqref{rotkur} during the simulation.  Panels E-H provide analogous information for a time shift of $\Delta t = -8$ hours.  Desynchronization is clearly observed during the application of each optimal stimulus in the raster plots from panels F and G.  Indeed, during both simulations, $r_{\rm rot}$ takes paths through the center in the complex plane.  }
\label{kuramotorecovery}
\end{figure}

Finally, the individual cell dynamics during recovery are considered in Figure 8.  For a time shift $\Delta t = +10$ hours that occurs instantaneously starting at $t = 0$ hours, traces of $E_i(t)$  for 30 representative oscillators are shown as colored lines in addition to $F(t)$ (solid black line) and a reference $F(t)$ that is fully entrained to the +10 hour time shifted light-dark cycle.  Each dot of the raster plot in Panel B represents the moment that a single cell in the simulation crosses the threshold $E_i(t) = 0.035$, and panel B shows the total light applied during the simulation.  The $\Delta L(t)$  is only applied during $t \in [0,48]$, after which, the 24 hour light-dark cycle remains.  Panel B clearly illustrates significant desynchronization that occurs during the application of the optimal stimulus.  After the stimulation is applied, the cells nearly reach their nominal level of synchronization with the desired phase.  This recovery is visualized in panel D, where $r_{\rm rot}$ is the Kuramoto order parameter \cite{kura84} taken in a rotating reference frame:
\begin{equation} \label{rotkur}
r_{\rm rot}(t) = \frac{1}{N} \sum_{k = 1}^N e^{i(\theta_k(t) -  2\pi t/24)},
\end{equation}
where $\theta_i(t)$ is the phase of oscillator $i$ at time $t$.  Here, $r_{\rm rot}$ will be referred to as the rotating order parameter.  Intuitively, Equation \eqref{rotkur} contains similar information to the standard definition of the order parameter; $|r_{\rm rot}|$ gives a sense of the cohesion of the oscillators, where $|r_{\rm rot}| = 1$ when the network is perfectly synchronized and $|r_{\rm rot}| = 0$ when the neurons are in a splay state.  The term $2\pi t/24$ is subtracted from the phase in \eqref{rotkur} so that when a time shift of $\Delta t$ is applied to the model, the argument of the complex number $r_{\rm rot}(t)$ shifts by $2 \pi \Delta t/24$ in steady state (i.e.,~once the model is fully reentrained to the shifted 24-hour light-dark cycle). Because information about the phases is not directly available, it is assumed that each oscillator reaches a phase of $\theta = 0$ when its value of $E_i$ crosses 0.035 with a positive slope, and all other phase values are linearly interpolated for times in between.   Panels E-H show similar information for a time shift of $\Delta t = -8$ hours when the optimal control obtained from the adaptive reduction is applied.

\section{Discussion and Conclusion} \label{concsec}
In this work, an optimal control strategy is developed for manipulating oscillation timing in the presence of a large magnitude entraining stimulus.  In contrast to other recently proposed strategies that leverage phase reduced models for control of oscillation timing \cite{moeh06}, \cite{nabi132}, \cite{mong19b}, \cite{wils17isored} the framework considered here is applicable for arbitrary large magnitude inputs.  Consequently, larger shifts in phase can be considered without sacrificing accuracy.  While this work is primarily motivated by the search for control protocols that hasten recovery from jet-lag in order to mitigate its negative side effects, the control formulation is quite general and could be implemented in other applications where rapid modification of oscillation timing is required.

In order to accurately account for the influence of large magnitude inputs in a reduced order setting, the recently developed phase-amplitude reduction framework is employed here.  By considering a family of limit cycles that emerge for different parameter sets, the nominal system parameters can be adaptively updated in order to limit the deviation from the nominal limit cycle, thereby minimizing error in the reduced order equations.  In general, the dimension of the adaptive reduction grows with both the number of isostable coordinates required and the number of the adaptive parameters.   However, as shown in this work, provided the amplitude corrections can be captured by a single isostable coordinate and $Q(\theta,p)$ is nonzero for all allowable $\theta$ and $p$, the adaptive parameter can be updated in such a way that the isostable coordinate dynamics are eliminated, yielding a reduced order set of equations with only two dimensions.  A subsequent optimal control formulation for manipulating the timing of oscillations based on Pontryagin's  minimum principle \cite{kirk98} is considered  in order to make comparisons with previously published results, however, a wide variety of optimal control formulations would be possible given the resulting two-dimensional reduced order representation of the dynamics.

For both the simple three-dimensional oscillator model \eqref{circmodel} and the model comprised of 3000 heterogeneous, noisy coupled oscillators \eqref{popmodel} the optimal control formulation using the adaptive reduction strategy vastly outperforms a strategy that uses a phase-only reduction \cite{moeh06} and a strategy that uses a first order accurate phase-amplitude reduction \cite{mong19b}.  When no entraining stimulus is considered (as is the case in results presented in Figure \ref{noconstraintresult}), either reduction framework is sufficient for identifying stimuli that can shift the oscillation timing.  However, the phase-only and first order phase-amplitude reduction strategy begins to degrade rapidly as the magnitude of the mandated shift grows.  Such degradation is not observed when the adaptive reduction is used, even when considering the largest possible phase shifts of $\Delta t = \pm 12$ hours (i.e.,~$\Delta \theta \approx \pm \pi$).  The incorporation of a large magnitude input that simulates the influence of a nominal  24-hour external light-dark cycle in the optimal control formulation renders the results when phase-only and first order accurate phase-amplitude reduced equations unusable; resulting recovery times in Figure \ref{constrainedresult} are similar to the uncontroled recovery for these formulations.  Conversely, optimal control waveforms obtained when using the adaptive reduction perform as desired with negligible errors.  A similar story emerges in the results from Figure \ref{recoverypopulation} for the coupled population model \eqref{popmodel}.  

A secondary focus of this work considers the identification of the necessary terms of an adaptive phase reduction using data-driven techniques.  The direct method \cite{izhi07}, \cite{neto12} is a well-established, data-driven technique for inferring phase reduced equations of the form \eqref{basicphase} when the underlying dynamical equations are unknown (e.g.,~in experimental systems).  While related methods have been proposed to infer numerical equations that characterize the behavior of the amplitude coordinates \cite{wils19phase}, \cite{wils20ddred}, there is still much work to be done.  As illustrated for the coupled population of limit cycle oscillators \eqref{popmodel}, provided that the adaptive system parameter is chosen appropriately, all of the terms of the adaptive phase reduction can be determined by applying a series of pulse and step function inputs and inferring the resulting change in the phase and isostable coordinates.  While this process requires more data than simply applying the direct method to compute the phase reduced equations \eqref{basicphase}, ignoring the extra information afforded by the adaptive reduction (e.g.,~by only implementing the first order accurate phase-amplitude reduction in the results shown in Figure \ref{recoverypopulation}) renders the resulting optimal control outputs useless.

When considering the problem of recovery from circadian misalignment, an interesting result observed here is that for both large magnitude time advances and time delays, the optimal control waveforms first desynchronize the individual oscillators within the population, effectively destroying the collective oscillation, before resynchronizing with the target phase of oscillation.  Such a notion is related to the observation from {\it in vito} experiments that the cells of the suprachiasmatic nucleus reentrain faster to shifts in environmental time if they are first desynchronized prior to the time shift \cite{an13}.  While the level of desynchronization was not explicitly considered by the optimal control formulation in Section \ref{popsec}, the resulting optimal control inputs significantly desynchronize the population during the phase shift.  In future work, it would be of interest to investigate the mathematical mechanisms by which desynchronization allows for more rapid shifts in the population phase of oscillation and whether these observations could be leveraged to develop better treatments for rapid recovery from circadian misalignment.  

While the adaptive phase reduction used here shows promise for control of oscillation timing in applications where large magnitude inputs are required, there are many opportunities for improvement.  From a theoretical perspective, little work has been done concerning the design of function $G_p$ of the adaptive reduction \eqref{fulladaptive} that is required to limit the magnitude of isostable coordinates $\psi_1,\dots,\psi_\beta$.  Of course, when only a small number of isostable coordinates are considered, as is the case in this work, adequate choices for $G_p$ are relatively easy to identify.  The same is not true when many isostable coordinates are considered and more consideration about the design of $G_p$ when multiple isostable coordinates are required would be warranted.  As a matter of practical implementation, it would be useful to identify strategies to experimentally infer the necessary terms of the adaptive reduction from \eqref{fulladaptive} that allow for arbitrary adaptive parameters to be considered.  A significant difficulty in inferring the terms of the adaptive reduction using data-driven methods is in the determination of $U_e(t,p,x)$.  Note that the definition from Equation \eqref{uedef}, $U_e$ requires knowledge of the underlying dynamical equations which are usually unknown in experimental applications.  The explicit assumption used in this work that $U(t) = \delta u(t)$, and   $F(x,p) = F(x) + p \delta$ where $\delta \in \mathbb{R}^N$ (i.e., that a shift in the nominal parameter $p$ provides the same effect as a shift in $u(t)$) provides a workaround to this limitation that allows $U_e$ to be computed exactly.  Nevertheless, this assumption severely limits the choice of the adaptive parameters and it would be useful to identify strategies that can infer the terms of the adaptive reduction for more general choices of the adaptive parameter.  Finally,  the amount of data required to infer the terms of the adaptive reduction using the methodology suggested in Section \ref{procsec} may be particularly onerous in some applications -- each parameter set considered requires a sufficient number of trials performed using both pulse and step function inputs to provide adequate data to fit functions to $z$, $i$, $D$, and $Q$ as illustrated in Figure \ref{calculateresponseterms}.  Strategies that can be used to infer the terms of the adaptive reduction more efficiently, such as those that use deep learning approaches \cite{lusc18}, \cite{yeun19} would certainly be of interest and will be investigated in future work.

This material is based upon the work supported by the National Science Foundation (NSF) under Grant No.~CMMI-1933583.


\end{document}